\documentclass{amsart}
\newcommand{\tp}[3]{\{#1#2#3\}}
   \newcommand{\tpc}[3]{\{#1,#2,#3\}}
\newcommand{\jbwst}{$JBW^*$-triple}
\newcommand{\csa}{$C^*$-algebra}

\theoremstyle{remark}
\newcommand{\tr}{\mbox{tr}\,}

\newcommand{\ip}[2]{\mbox{$(#1|#2)$}}
\newcommand{\CC}{{\bf C}}
\newcommand{\NN}{\bf N}

\newtheorem{lemma}{Lemma}[section]
\newtheorem{theorem}{Theorem}
\newcommand{\pf}{{\it Proof}.}
\newtheorem{proposition}[lemma]{Proposition}
\newtheorem{remark}[lemma]{Remark}
\newtheorem{corollary}[lemma]{Corollary}

\newtheorem{example}{Example}
\newtheorem{defn}{Definition}

\newcommand{\jbst}{$JB^*$-triple}
\newcommand{\jcst}{$JC^*$-triple}
\newcommand{\jwst}{$JW^*$-triple}

\begin{document}
\title[
Contractively complemented Hilbertian operator
spaces]{Classification of contractively complemented Hilbertian
operator spaces}
\author{Matthew Neal}
\address{Department of Mathematics, Denison University, Granville,
Ohio 43023}
 \email{nealm@denison.edu}

\author{Bernard Russo}
\address{Department of Mathematics, University of California,
Irvine, California 92697-3875}
 \email{brusso@math.uci.edu}


\date{August 26, 2005}



\begin{abstract}
We construct some separable infinite dimensional homogeneous
Hilbertian operator spaces $H_\infty^{m,R}$ and $H_\infty^{m,L}$,
which generalize the row and column spaces $R$ and $C$ (the case
$m=0$). We show that separable infinite-dimensional Hilbertian
$JC^*$-triples are completely isometric to an element of the set of
(infinite) intersections of these spaces . This set includes the
operator spaces $R$, $C$, $R \cap C$, and the space $\Phi$ spanned
by creation operators on the full anti-symmetric Fock space.  In
fact, we show that $H_{\infty}^{m,R}$ (resp. $H_\infty^{m,L}$) is
completely isometric to the space of creation (resp. annihilation)
operators on the $m$ (resp. $m+1$) anti-symmetric tensors of the
Hilbert space. Together with the finite-dimensional case studied in
\cite{NeaRus05}, this gives a full operator space classification of
all rank-one $JC^*$-triples in terms of creation and annihilation
operator spaces.

We use the above to show that all contractive projections on a
C*-algebra $A$ with infinite dimensional Hilbertian range are
``expansions'' (which we define precisely) of normal contractive
projections from $A^{\ast\ast}$ onto a Hilbertian space which is
completely isometric to $R$, $C$, $R \cap C$, or $\Phi$.
This generalizes the well known result, first proved for $B(H)$ by
Robertson in \cite{R}, that all Hilbertian operator spaces that are
{\it completely} contractively complemented in a C*-algebra are
completely isometric to $R$ or $C$.  We also compute various
completely bounded Banach-Mazur distances between these spaces, or
$\Phi$.
\end{abstract}

\maketitle

\section{Preliminaries}

The goals of the present paper are to classify all infinite
dimensional rank 1 JC*-triples up to complete isometry (Theorem
\ref{thm:1} in section 2) and then use that result to give a
suitable ``classification'' of all Hilbertian operator spaces which
are contractively complemented in a C*-algebra or normally
contractively complemented in a W*-algebra (Theorems 2 and 3 in
section 3). In particular, we show  that these spaces are
``essentially'' $R$, $C$, $R\cap C$, or $\Phi$ modulo a
``degenerate'' piece. 

In Theorem 4 in section 4 we compute various completely bounded
Banach-Mazur distances between these JC*-triples.  In Theorem 5 in
section 5, we show that all of these JC*-triples in the separable
infinite dimensional and finite dimensional cases can be represented
completely isometrically as creation and annihilation operator
spaces on pieces of the anti-symmetric Fock space.

In the rest of this section, we give some background on operator
space theory and on $JC^*$-triples.

\subsection{Operator spaces} Operator space theory is a
non-commutative or quantized theory of Banach spaces. By definition,
an operator space is a Banach space together with an isometric
linear embedding into $B(H)$, the bounded linear operators on a
complex Hilbert space. While the objects are obviously the Banach
spaces themselves, the more interesting aspects concern the
morphisms, namely, the completely bounded maps. These are defined by
considering an operator space as a subspace $X$ of $B(H)$.
 Its {\it
operator space structure} is then given by the sequence of norms on
the set of matrices $M_n(X)$ with entries from $X$, determined by
the identification $M_n(X) \subset M_n(B(H))=B(H\oplus H \oplus
\cdots \oplus H)$. A linear mapping $\varphi:X\rightarrow Y$ between
two operator spaces is {\it completely bounded} if the induced
mappings $\varphi_n:M_n(X)\rightarrow M_n(Y)$ defined by
$\varphi_n([x_{ij}])=[\varphi(x_{ij})]$ satisfy
$\|\varphi\|_{\mbox{cb}}:=\sup_n\|\varphi_n\|<\infty$.

 Operator space
theory has its origins in the work of Stinespring in the 1950s, and
Arveson in the 1960s.  Many tools were developed in the 1970s and
1980s by a number of operator algebraists, and an abstract framework
was developed in 1988 in the thesis of Ruan. All definitions,
notation, and results used in this paper can be found in recent
accounts of the subject, namely (in chronological order)
\cite{EffRau00},\cite{Paulsen02},\cite{Pisier03},\cite{BleLeM04}.
Let us just recall that a completely bounded map is a {\it complete
isomorphism} if its inverse exists and is completely bounded. Two
operator spaces are {\it completely isometric} if there is a linear
isomorphism $T$ between them with $\|T\|_{\mbox{cb}}=
\|T^{-1}\|_{\mbox{cb}}=1$. We call $T$ a {\it complete isometry} in
this case. Other important types of morphisms in this category are
complete contractions ($\|\varphi\|_{\mbox{cb}}\le 1$) and complete
semi-isometries (:= isometric complete contraction).

Examples of completely bounded maps  are the restriction to a
subspace of a $C^*$-algebra of a *-homomorphism and multiplication
by an fixed element. It is a fact that every completely bounded map
is essentially a product of these two examples, \cite[Th.
1.6]{Pisier03}. The space $CB(X,Y)$ of completely bounded maps
between operator spaces $X$ and $Y$ is a Banach space with the
completely bounded norm $\|\cdot \|_{\mbox{cb}}$.

Analogous to the Banach-Mazur distance for Banach spaces, the class
of all operator spaces can be made into a metric space by using the
logarithm of the {\it completely bounded} Banach-Mazur distance:
\[
\mbox{d}_{\mbox{cb}}(E,F)=\inf\{\|u\|_{\mbox{cb}}\cdot
\|u^{-1}\|_{\mbox{cb}}\ ;\  u:E\rightarrow F\mbox{ complete
isomorphism} \}.
\]

Two important examples of Hilbertian operator spaces (:= operator
spaces isometric to Hilbert space) are the row and column spaces
 $R,\ C$, and their finite-dimensional versions $R_n,\ C_n$.
These are defined as follows. In the matrix representation for
$B(\ell_2)$,  {\it column Hilbert space}
$C:=\overline{\mbox{sp}}\{e_{i1}:i\ge 1\}$ and {\it row Hilbert
space} $R:=\overline{\mbox{sp}}\{e_{1j}:j\ge 1\}$. Their finite
dimensional versions are $C_n=\mbox{sp}\{e_{i1}:1\le i\le n\}$ and
$R_n=\mbox{sp}\{e_{1j}:1\le j\le n\}$. Here of course $e_{ij}$ is
the operator defined by the matrix with a 1 in the $(i,j)$-entry and
zeros elsewhere. Although $R$ and $C$ are Banach isometric, they are
not completely isomorphic ($\mbox{d}_{\mbox{cb}}(R,C)= \infty$); and
$R_n$ and $C_n$, while completely isomorphic,
 are not completely isometric.
In fact, it is known that $\mbox{d}_{\mbox{cb}}(R_n,C_n)= n$.

$R,\ C,\ R_n,\ C_n$ are examples of {\it homogeneous} operator
spaces, that is, operator spaces $E$ for which $\forall
u:E\rightarrow E$, $\|u\|_{\mbox{cb}}=\|u\|$. Another important
example of an Hilbertian homogeneous operator space is $\Phi(I)$.
The  space $\Phi(I)$ is defined by
$\Phi(I)=\overline{\mbox{sp}}\{V_i:i\in I\}$, where the $V_i$ are
bounded operators on a Hilbert space satisfying the canonical
anti-commutation relations. In some special cases, the notations
 $\Phi_n:=\Phi(\{1,2,\ldots,n\})$, and
$\Phi=\Phi(\{1,2,\ldots \})$ are used. For more properties of this
space and related constructs, see \cite[9.3]{Pisier03}.

Two more examples of homogeneous operator spaces are $\min(E)$,
$\max(E)$, where $E$ is any Banach space. For any such $E$, the
operator space structure of $\min(E)$ is defined by the embedding of
$E$ into the continuous functions on the unit ball of $E^*$ in the
weak*-topology, namely, $\|(a_{ij})\|_{M_n(\min(E))}=\sup_{\xi\in
B_{E^*}}\|(\xi(a_{ij}))\|_{M_n}$. The operator space structure of
 $\max(E)$ is given by
$$\|(a_{ij})\|_{M_n(\max(E))}=\sup\{\|(u(a_{ij}))\|_{M_n(B(H_u))}:u:E\rightarrow
B(H_u),\ \|u\|\le 1 \}.$$ More generally, if $F$ and $G$ are
operator spaces, then in $F\stackrel{u}{\longrightarrow}\min(E)$,
$\|u\|_{\mbox{cb}}=\|u\|$, and in $
\max(E)\stackrel{v}{\longrightarrow}G$, $\|v\|_{\mbox{cb}}=\|v\|$.
The notations $\min(E)$ and $\max (E)$ are justified by the fact
that for any Banach space $E$, the identity map on $E$ is completely
contractive in
 $\max(E)\rightarrow
 E\rightarrow
\min(E)$.

By analogy with the classical Banach spaces $\ell_p,\ c_0,\ L_p,\
C(K)$ (as well as their ``second generation'', Orlicz, Sobolev,
Hardy, Disc algebra, Schatten $p$-classes), we can consider the
 (Hilbertian) operator spaces $R,\ C,\ \min(\ell_2),\
\max(\ell_2),\ OH,\ \Phi$, as well as their finite dimensional
versions
 $R_n,\ C_n,\ \min(\ell_2^n),\ \max(\ell_2^n),\ OH_n,\ \Phi_n$, as
``classical operator
 spaces''. Among these spaces, only the spaces $R,C$, and $\Phi$ play
important roles in
 this paper. (For the definition and properties of the space called $OH$,
 see \cite[Chapter 7]{Pisier03}.)
The classical operator spaces are mutually completely
non-isomorphic. If $E_n,F_n$ are $n$-dimensional versions, then
$\mbox{d}_{\mbox{cb}}(E_n,F_n)\rightarrow\infty$, \cite[Ch.\
10]{Pisier03}.

We propose to add to this list of classical operator spaces the
Hilbertian operator spaces
 $H_\infty^{m,R}$ and
 $H_\infty^{m,L}$ constructed here, as well as their
 finite-dimensional versions $H_n^k$  studied in \cite{NeaRus03} and
 \cite{NeaRus05}.
Like the space $\Phi$, the spaces $H_\infty^{m,R}$, $H_\infty^{m,L}$
and $H_n^k$ can be
 represented up to complete isometry as spaces of creation
 operators or annihilation operators on  anti-symmetric Fock
 spaces (\cite[Lemma 2.1]{NeaRus05} and Theorem~\ref{thm:5} below).

Let us recall from \cite[Sections 6,7]{NeaRus03} the construction of
the spaces $H_n^k$, $1\le k\le n$. Let $I$ denote a subset of
$\{1,2,\ldots,n\}$ of cardinality $|I|=k-1$. The number of such $I$
is $q:={n\choose k-1}$. Let $J$ denote a subset of
$\{1,2,\ldots,n\}$ of cardinality $|J|=n-k$. The number of such $J$
is $p:={n\choose n-k}$. We  assume that each
$I=\{i_1,\ldots,i_{k-1}\}$ is such that $i_1<\cdots <i_{k-1}$, and
that if $J=\{j_1,\ldots,j_{n-k}\}$, then $j_1<\cdots<j_{n-k}$.

The space $H_n^k$ is the linear span of matrices $b_i^{n,k}$, $1\le
i\le n$, given by
\[
b_i^{n,k}=\sum_{I\cap J=\emptyset,(I\cup
J)^c=\{i\}}\epsilon(I,i,J)e_{J,I},
\]
where $e_{J,I}=e_J\otimes e_I=e_Je_I^t\in
M_{p,q}(\CC)=B(\CC^q,\CC^p)$, and $\epsilon(I,i,J)$  is the
signature of the permutation taking
$(i_1,\ldots,i_{k-1},i,j_1,\ldots,j_{n-k})$ to $(1,\ldots,n)$. Since
the $b_i^{n,k}$ are the image under a triple isomorphism (actually
ternary isomorphism) of a rectangular grid in a \jwst\ of rank one,
they form an orthonormal basis for $H_n^k$ (cf. \cite[subsection 5.3
and section 7]{NeaRus03}).

The following definition from \cite[2.7]{Pisier03} plays a key role
in this paper. If $E_0\subset B(H_0)$ and $E_1\subset B(H_1)$ are
operator spaces whose underlying Banach spaces form a compatible
pair in the sense of interpolation theory, then the Banach space
$E_0\cap E_1$ (with the norm $\|x\|_{E_0\cap E_1}=\max
(\|x\|_{E_0},\|x\|_{E_1})$) equipped with the operator space
structure given by the embedding $E_0\cap E_1 \ni x\mapsto (x,x)\in
E_0\oplus E_1\subset B(H_0\oplus H_1)$ is called the {\it
intersection} of $E_0$ and $E_1$ and is denoted by $E_0\cap E_1$. We
note, for examples, that $\cap_{k=1}^n H_n^k=\Phi_n$
(\cite{NeaRus05}) and the space $R\cap C$ is defined relative to the
embedding of $C$ into itself and $R$ into $C$ given by the transpose
map (\cite[p.\ 184]{Pisier03}). The definition of intersection
extends easily to arbitrary families of compatible operator spaces
(cf. Theorem~\ref{thm:1} below).

\begin{lemma}\label{lem:2.1}
Let $H$ be an Hilbertian operator space, and suppose that every
finite dimensional subspace of $H$ is homogeneous.  Then $H$ itself
is homogeneous.
\end{lemma}
\pf\ Let $\phi$ be any unitary operator on $H$.  According to the
first statement of \cite[Prop.9.2.1]{Pisier03}, it suffices to prove
that $\phi$ is a complete isometry.

Let $F$ be any finite dimensional subspace of $H$ and let $G$ be the
subspace spanned by $F\cup \phi(F)$.  By the second statement of
\cite[Prop.9.2.1]{Pisier03}, $F$ and $\phi(F)$, being of the same
dimension as subspaces of the homogeneous space $G$, are completely
isometric, and  $\phi|F$ is a complete isometry.

Now let $[x_{ij}]\in M_n(H)$.  Then $\{x_{ij},\phi(x_{ij}):1\le
i,j\le n\}$ spans a finite dimensional subspace $F$ of $H$, and
\[
\|\phi_n([x_{ij}])\|_{M_n(H)}=\|\phi_n([x_{ij}])\|_{M_n(F)}=
\|[x_{ij}]\|_{M_n(F)}=\|[x_{ij}]\|_{M_n(H)}.\qed
\]

\subsection{Rank one \jcst s}

A \jcst\ is a norm closed complex linear subspace  of $B(H,K)$
(equivalently, of a $C^*$-algebra) which is closed under the
operation $a\mapsto aa^*a$. \jcst s were defined and studied (using
the name $J^*$-algebra)
 as a generalization of \csa s by Harris \cite{Harris73} in
connection with function theory on infinite dimensional bounded
symmetric domains. By a polarization identity, any \jcst\ is closed
under the triple product
\begin{equation}\label{eq:product}
(a,b,c)\mapsto \tp{a}{b}{c}:=\frac{1}{2}(ab^*c+cb^*a),
\end{equation}
under which it becomes a Jordan triple system.   A linear map which
preserves the triple product (\ref{eq:product}) will be called a
{\it triple homomorphism}. Cartan factors are examples of \jcst s,
as are \csa s, and Jordan \csa s. Cartan factors are defined for
example in \cite[Section 1]{NeaRus03}. We shall only make use of
Cartan factors of type 1, that is, spaces of the form $B(H,K)$ where
$H$ and $K$ are complex Hilbert spaces.

A special case of a \jcst\ is a {\it ternary algebra}, that is, a
subspace of $B(H,K)$ closed under the {\it ternary product}
$(a,b,c)\mapsto ab^*c$. A {\it ternary homomorphism} is a linear map
$\phi$ satisfying $\phi(ab^*c)=\phi(a)\phi(b)^*\phi(c)$. These
spaces are also called ternary rings of operators and abbreviated
TRO. They have been studied both concretely in \cite{Hestenes62} and
abstractly in \cite{Zettl83}. Given a TRO $M$, its left (resp.
right) linking \csa\ is defined to be the norm closed span of the
elements $ab^*$ (resp. $a^*b$) with $a,b\in M$. Ternary isomorphic
TROs have isomorphic left and right linking algebras.

TROs have come to play a key role in operator space theory, serving
as the algebraic model in the category. Recall that the algebraic
models for the categories of order-unit spaces, operator systems,
and Banach spaces, are respectively Jordan $C^*$-algebras, \csa s,
and \jbst s.  Indeed, for TROs, a ternary isomorphism is the same as
a complete isometry.

If $v$ is a partial isometry in a \jcst\  $M\subset B(H,K)$, then
the projections $l=vv^*\in B(K)$ and $r=v^*v\in B(H)$ give rise to
(Peirce) projections $P_k(v):M\rightarrow M,\ k=2,1,0$  as follows;
for $x\in M$,
\[
P_2(v)x=lxr\quad , \quad P_1(v)x=lx(1-r)+(1-l)xr\quad , \quad
P_0(v)x=(1-l)x(1-r).
\]
The projections $P_k(v)$ are contractive,  and their ranges, called
Peirce spaces and denoted by $M_k(v)$, are $JC^*$-subtriples of $M$
satisfying $M=M_2(v)\oplus M_1(v)\oplus M_0(v)$.

A partial isometry $v$ is said to be {\it minimal} in $M$ if
$M_2(v)={\CC} v$. This is equivalent to $v$ not being the sum of two
non-zero orthogonal partial isometries. Recall that two partial
isometries $v$ and $w$ (or any two Hilbert space operators) are
orthogonal if $v^*w=vw^*=0$. Orthogonality of partial isometries $v$
and $w$ is equivalent to $v\in M_0(w)$ and will be denoted by
$v\perp w$. Each finite dimensional \jcst\ is the linear span of its
minimal partial isometries. More generally, a \jcst\ is defined to
be {\it atomic} if it is the weak closure of the span of its minimal
partial isometries. In this case, it has a predual and is called a
{\it \jwst}. The {\it rank} of a \jcst\ is the maximum number of
mutually orthogonal minimal partial isometries. For example, the
rank of the Cartan factor $B(H,K)$ of type 1  is the minimum of the
dimensions of $H$ and $K$; and the rank of the Cartan factor of type
4 (spin factor) is 2.

In a JC*-triple, there is a natural ordering on partial isometries.
We write $v \le w$ if $vw^{*}v=v$; this is equivalent to $vv^*\le
ww^*$ and $v^*v\le w^*w$. Moreover, if $v\le w$, then there exists a
partial isometry $v'$ orthogonal to $v$ with $w=v+v'$.

Another relation between two partial isometries that we shall need
is defined in terms of the Peirce spaces as follows. Two partial
isometries $v$ and $w$ are said to be {\it collinear} if $v\in
M_1(w)$ and $w\in M_1(v)$, notation $v\top w$. Let $u,v,w$ be
partial isometries. The following is part of \cite[Lemma
5.4]{NeaRus03}, and is referred to as {\it ``hopping''}: If $v$ and
$w$ are each collinear with $u$, then $uu^*vw^*=vw^*uu^*$ and
$u^*uv^*w=v^*wu^*u$.  If $u,v,w$ are mutually collinear partial
isometries, then $\tp{u}{v}{w}=0$.

\jcst s of arbitrary dimension occur naturally in functional
analysis and in holomorphy. A
 special case of a theorem of Friedman and Russo \cite[Theorem 2]{FriRus85}
 states that if $P$ is a contractive
projection on a $C^*$-algebra $A$, then there  is a linear isometry
of the range $P(A)$ of $P$ onto a $JC^*$-subtriple of $A^{**}$. A
special case of a theorem of Kaup \cite{Kaup83} gives a bijective
correspondence between Cartan factors and irreducible bounded
symmetric domains in complex Banach spaces.

Contractive projections play a ubiquitous role in the structure
theory of the abstract analog of \jcst s (called $JB^*$-triples). Of
use to us will be both of the following two conditional expectation
formulas for a contractive projection $P$ on a JC*-triple $M$ (which
are valid for \jbst s) (\cite[Corollary
1]{FriRus84}):
\begin{equation}\label{eq:CE}
P\tpc{Px}{Py}{Pz}=P\tpc{Px}{Py}{z}=P\tpc{Px}{y}{Pz},\quad (x,y,z\in
M).
\end{equation}

By a special case of  \cite[Cor.,p.308]{DanFri87}, every \jwst\ of
rank one is isometric to a Hilbert space and every maximal collinear
family of partial isometries corresponds to an orthonormal basis.
Conversely, every Hilbert space with the abstract triple product
$\tp{x}{y}{z}:=(\ip{x}{y}z+\ip{z}{y}x)/2$ can be realized as a
\jcst\ of rank one in which every orthonormal basis forms a maximal
family of mutually collinear minimal partial isometries.

\section{Operator space structure of Hilbertian JC*-triples}\label{sec:coor}

\subsection{Hilbertian JC*-triples: The spaces $H_\infty^{m,R}$
and $H_\infty^{m,L}$}\label{sec:2.1}

The general setting for the next two sections will be the following:
$Y$ is a $JC^*$-subtriple of $B(H)$ which is Hilbertian in the
operator space structure arising from $B(H)$, and
$\{u_i:i\in\Omega\}$ is an orthonormal basis consisting of a maximal
family of mutually collinear partial isometries of $Y$. Note that
the $u_i$ are each minimal in $Y$, but not necessarily minimal in
any $JC^*$-triple containing $Y$.

We let $T$ and $A$ denote the TRO and the $C^*$-algebra respectively
generated by $Y$. For any subset $G\subset\Omega$,
$(uu^*)_G:=\prod_{i\in G}u_iu_i^*$ and $(u^*u)_G:=\prod_{i\in
G}u_i^*u_i$. The elements $(uu^*)_G$ and $(u^*u)_G$ lie in the weak
closure of $A$ and more generally in  the left and right linking von
Neumann algebras of $T$.

In the following lemma, parts (a) and (a$^\prime$) justify the
definitions of the integers $m_R$ and $m_L$ in parts (b) and
(b$^\prime$).

\begin{lemma}\label{lem:3.1}
Let $Y$ be an Hilbertian operator space which is a $JC^*$-subtriple
of $B(H)$ and let $\{u_i:i\in\Omega\}$ be an orthonormal basis
consisting of a maximal family of mutually collinear partial
isometries of $Y$.
\begin{description}
\item[(a)] If
$(uu^*)_{\Omega-F}=0$ for some finite set $F\subset \Omega$, then,
$(uu^*)_{\Omega-G}=0$ for every finite set $G$ with the same
cardinality as $F$.
\item[(a$^\prime$)] If
$(u^*u)_{\Omega-F}=0$ for some finite set $F\subset \Omega$, then,
$(u^*u)_{\Omega-G}=0$ for every finite set $G$ with the same
cardinality as $F$.
\item[(b)] Assume $(uu^*)_{\Omega-F}\ne 0$ for some finite set $F$. Let
$m_R$ be the smallest nonnegative integer with $(uu^*)_{\Omega-F}\ne
0$ for every $F$ with cardinality $m_R$. Define
$p_R=\sum_{|F|=m_R}(uu^*)_{\Omega-F}$. Then the maps $y\mapsto p_Ry$
and $y\mapsto (1-p_R)y$ are completely contractive triple
isomorphisms of $Y$ onto rank one subtriples of the weak closure of
$T$ in $B(H)$. Moreover, $p_RY\perp (1-p_R)Y$.
\item[(b$^\prime$)] Assume $(u^*u)_{\Omega-F}\ne 0$ for some finite set
$F$. Let $m_L$ be the smallest nonnegative integer with
$(u^*u)_{\Omega-F}\ne 0$ for every $F$ with cardinality $m_L$.
Define $p_L=\sum_{|F|=m_L}(u^*u)_{\Omega-F}$. Then the maps
$y\mapsto yp_L$ and $y\mapsto y(1-p_L)$ are completely contractive
triple isomorphisms of $Y$ onto rank one subtriples of the weak
closure of $T$ in $B(H)$. Moreover, $Yp_L\perp Y(1-p_L)$.
\item[(c)] In case (b), let $w_i=p_Ru_i$ and let
$m'_R$ be the smallest nonnegative integer with
$(ww^*)_{\Omega-F}\ne 0$ for all $F$ with cardinality $m'_R$. Then
$m_R'$ exists, and $m'_R=m_R$. Furthermore, $(w^*w)_G\ne 0$ for
$|G|=m_R+1$ and $(w^*w)_G= 0$ for $|G|=m_R+2$. Thus, if we define
$k_R$ to be the smallest integer $k$ such that $(w^*w)_G\ne 0$ for
$|G|=k$, then $k_R=m_R+1$.
\item[(c$^\prime$)] In case (b$^\prime$), let $w_i=u_ip_L$ and let
$m'_L$ be the smallest nonnegative integer with
$(w^*w)_{\Omega-F}\ne 0$ for all $F$ with cardinality $m'_L$. Then
$m_L'$ exists, and $m'_L=m_L$. Furthermore, $(ww^*)_G\ne 0$ for
$|G|=m_L+1$ and $(ww^*)_G= 0$ for $|G|=m_L+2$.   Thus, if we define
$k_L$ to be the smallest integer $k$ such that $(ww^*)_G\ne 0$ for
$|G|=k$, then $k_L=m_L+1$.
\end{description}
\end{lemma}
\pf\ The proofs of (a) and (a$^\prime$) are identical to the proof
given in \cite[Lemma 5.8]{NeaRus03}. The fact that the set
$\Omega-F$ is infinite has no effect on the proof in
\cite{NeaRus03}.

The proofs of (b) and (b$^\prime$) are identical to the proof given
in \cite[Lemma 5.9]{NeaRus03}. The facts that the set $\Omega-F$ is
infinite and that the sums defining the projections $p_R$ and $p_L$
are infinite have no effect on the proof in \cite{NeaRus03}.

We now prove (c), the proof of (c$^\prime$) being entirely similar.
For any finite set $F\subset\Omega$,
\begin{eqnarray*}
(ww^*)_{\Omega-F}&=&\prod_{i\in\Omega-F}\left(\sum_{|G|=m_R}(uu^*)_{\Omega-G}
\right)u_iu_i^*\left(\sum_{|H|=m_R}(uu^*)_{\Omega-H}\right)\\
&=&\prod_{i\in\Omega-F}\left(\sum_{|G|=m_R,i\in\Omega-G}(uu^*)_{\Omega-G}\right)=\sum_{G\subset
F,|G|=m_R} (uu^*)_{\Omega-G}.
\end{eqnarray*}
 From this it follows that $(ww^*)_{\Omega-F}=0$ if $|F|<m_R$ and
 that $(ww^*)_{\Omega-F}=(uu^*)_{\Omega-F}\ne 0$ if $|F|=m_R$.
 This proves that $m_R'=m_R$, that is
$(ww^*)_{\Omega-F}=0\Leftrightarrow |F|<m_R$.

Now let $|F|=r$ and for convenience, suppose that
$F=\{1,2,\ldots,r\}$.  Then
\[
(w^*w)_F=(w^*w)_{\{1,2,\ldots,r\}}=\sum
u_1^*(uu^*)_{\Omega-F_1}u_1u_2^*(uu^*)_{\Omega-F_2} u_2u_3^*\cdots
u_r^*(uu^*)_{\Omega-F_r}u_r,
\]
where the sum is over all $|F_j|=m_R,j\in\Omega-F_j,F-\{j\}\subset
F_j$ (by ``hopping''), and $j=1,2,\ldots r$. Every term in this sum
is zero if $r-1>m_R$, that is $r\ge m_R+2$.  Further, if $r=m_R+1$,
there is only one term, namely $x:=$
\begin{eqnarray*}
(w^*w)_{\{1,2,\ldots,m+1\}}&=&
u_1^*(uu^*)_{\Omega-\{2,3,\ldots,m+1\}}u_1u_2^*(uu^*)_{\Omega-\{1,3,4,\ldots,m+1\}}u_2u_3^*
\times\cdots\times\\
&&(uu^*)_{\Omega-\{1,2,\ldots,m-1,m+1\}}u_mu_{m+1}^*(uu^*)_{\Omega-\{1,2,\ldots,m\}}u_{m+1}.
\end{eqnarray*}
which by a sequence of ``hoppings'' becomes
\[
x=(u^*u)_{\{1,2,\ldots,m\}}u_{m+1}^*(uu^*)_{\Omega-\{1,2,\ldots,m+1\}}u_{m+1}.
\]
In turn, using the collinearity of the $u_k$, this becomes
\[
x=u_{m+1}^*(uu^*)_{\Omega-\{1,2,\ldots,m+1\}}u_{m+1}.
\]
Thus, if $x=0$, then
$0=u_{m+1}xu_{m+1}^*=(uu^*)_{\Omega-\{1,2,\ldots,m\}}$, a
contradiction. \qed

\medskip

 Our goal for the remainder of
this section is to give a completely isometric representation for
the spaces $p_RY$ and $Yp_L$ in parts (b) and (b$^\prime$) of Lemma
\ref{lem:3.1}.  This will be achieved via a coordinatization
procedure which we now describe.

In the following, let us restrict to the special case that $Y$ is a
Hilbertian JC*-triple which satisfies the properties of $p_RY$ in
Lemma \ref{lem:3.1} part (c). For notational convenience, let
$m=m_R$. Thus $(u^*u)_G\ne 0$ for $|G|\le m+1$ and $(u^*u)_G=0$ for
$|G|\ge m+2$.

Analogous to \cite[Def.\ 6.1]{NeaRus03}, we are going to define
elements which are indexed by an arbitrary pair of subsets $I,J$ of
$\Omega$ satisfying
\begin{equation}\label{eq:12}
|\Omega-I|=m+1,\quad |J|=m.
\end{equation}
Here and throughout the rest of this paper, $|F|$ denotes the
cardinality of the finite set $F$.

 The set $I\cap J$ is finite,
and if $|I\cap J|=s\ge 0$, then $|(I\cup J)^c|=s+1$. Let us write
$I\cap J=\{d_1,\ldots,d_s\}$ and
 $ (I\cup J)^c=\{c_1,\ldots,c_{s+1}\}$, and
let us agree (for the moment) that there is a natural linear
ordering on $\Omega$ such that $c_1<c_2<\cdots <c_{s+1}$ and
$d_1<d_2<\cdots<d_s$.

With the above notation, we define
\begin{equation}\label{eq:(I,J)}
u_{IJ}=u_{I,J}=(uu^*)_{I- J}u_{c_1}u_{d_1}^*u_{c_2}u_{d_2}^*\cdots
u_{c_s}u_{d_s}^*u_{c_{s+1}}(u^*u)_{J- I}.
\end{equation}

Note that in general $I-J$ is infinite and $J-I$ is finite so that
$u_{I,J}$ lies in the weak closure of $T$.

 In the special case of {\rm (\ref{eq:(I,J)})} where
 $I\cap J=\emptyset$, we have $s=0$ and $u_{I,J}$
has the form
\[
u_{I,J}=(uu^*)_Iu_c(u^*u)_J,
\]
where
 $I\cup J\cup \{c\}=\Omega$ is a partition of $\Omega$.
As in \cite{NeaRus03}, we call such an element a {\rm ``one''}, and
denote it also by $u_{I,c,J}$.

The proof of the following lemma, which is the analog of \cite[Lemma
6.6]{NeaRus03}, is complicated by the fact that the sets $I$ are
infinite if $\Omega$ is infinite.

\begin{lemma}\label{lem:3.0}
Let $Y$ be an Hilbertian operator space which is a $JC^*$-subtriple
of $B(H)$ and let $\{u_i:i\in\Omega\}$ be an orthonormal basis
consisting of a maximal family of mutually collinear partial
isometries of $Y$. Assume that $Y$ satisfies the properties of
$p_RY$ in Lemma \ref{lem:3.1} part (c) with $m=m_R$, that is,
$(u^*u)_G\ne 0$ for $|G|\le m+1$ and $(u^*u)_G=0$ for $|G|\ge m+2$.
For any $c\in \Omega$,
\begin{equation}\label{eq:sum}
u_c=\sum_{I,J}u_{I,J}=\sum_{I,J}u_{I,c,J}
\end{equation}
where the sum is taken over all disjoint $I,J$ satisfying {\rm
(\ref{eq:12})} and not containing $c$, and converges weakly in the
weak closure  of $T$.
\end{lemma}
\pf\ The proof of \cite[Prop.\ 6.3]{NeaRus03} remains valid in our
context insofar as  $\{u_{I,J}\}$ is a collection of pairwise
orthogonal partial isometries in the weak closure  of the ternary
envelope $T$ of $Y$. Since $u_c^*u_c$ commutes with $(u^*u)_J$,
$u_{I,J}^*u_{I,J}u_c^*u_c=u_{I,J}^*u_{I,J}$, so that
$u_{I,J}^*u_{I,J}\le u_c^*u_c$ and similarly $u_{I,J}u_{I,J}^*\le
u_cu_c^*$ so that $\sum u_{I,J}\le u_c$.   To prove (\ref{eq:sum}),
we proceed as follows.

Let us write $u_c=v_c+w_c$, where $w_c=\sum_{I,c,J} Iu_cJ$ and for
example, $Iu_cJ$ is shorthand for $u_{I,c,J}=(uu^*)_Iu_c(u^*u)_J$,
and $v_c$ is a partial isometry orthogonal to $w_c$. We shall show
that $v_c=0$. From the simple facts that $Iu_cJ=u_cJ$ and
$w_cJ=Iu_cJ$, it follows that $v_cJ=0$. Similarly, $Iv_c=0$. From
this, it follows that for each pair $i\ne j$, $v_i\in M_1(u_j)$.
Indeed,
\begin{eqnarray*}
v_i+w_i&=&u_i=u_iu_j^*u_j+u_ju_j^*u_i\\
&=&v_iu_j^*u_j +w_iu_j^*u_j+u_ju_j^*v_i+u_ju_j^*w_i\\
&=&v_iu^*_ju_j+\sum_{j\in J}u_iJ+u_ju_j^*v_i+\sum_{j\not\in
J}u_iJ\\
&=& v_iu^*_ju_j+u_ju_j^*v_i+w_i.
\end{eqnarray*}
Hence $v_i=v_iu^*_ju_j+u_ju_j^*v_i$, so $v_i\in M_1(u_j)$.

Next, we observe that $v_i$ is orthogonal to $w_j$ for every $i$ and
$j$.  For $i=j$ this is clear by definition. For $i\ne j$, we have
$v_i(Ju_j^*I)=(v_iJ)u_j^*I=0$ and similarly $Ju_j^*Iv_i=0$ so that
$v_iw_j^*=w_j^*v_i=0$.

It now follows that $v_i\top v_j$ for $i\ne j$. Indeed,
\[
v_j=v_ju_i^*u_i+u_iu_i^*v_j=v_jv_i^*v_i+v_jw_i^*w_i+v_iv_i^*v_j+w_iw_i^*v_j=
v_jv_i^*v_i+v_iv_i^*v_j.
\]

Let us adopt the notation $J_v$ for $(v^*v)_J=\prod_{j\in
J}v_j^*v_j$. (What we previously denoted by $J$ would now be denoted
by $J_u=(u^*u)_J$.)  We know that $v_iJ_v=v_iJ_u=0$. Suppose that
$v_iJ'_v\ne 0$ for some $J'\subset J$. We will show that
$J'=\emptyset$.  In the first place, $v_iJ'_v=I'_vv_iJ'_v$, since
letting $I'=\Omega-(J'\cup\{i\})=\{i_\alpha:\alpha\in \Lambda\}$
say,
\[
v_iJ'_v=v_i(v^*v)_{J'}=(v_{i_\alpha}v_{i_\alpha}^*v_i+v_iv_{i_\alpha}^*v_{i_\alpha})
(v^*v)_{J'}
=v_{i_\alpha}v_{i_\alpha}^*v_i(v^*v)_{J'}=\cdots=I'_vv_iJ'_v.
\]
In the second place, by the orthogonality of $v_j$ and $w_j$,
\[
I'_uu_j=(uu^*)_{I'}u_j=[\prod_{i\in
I'}(v_i+w_i)(v_i^*+w_i^*)](v_j+w_j)=I'_vv_j+I'_ww_j.
\]
Hence, $I'_vv_j=I'_uu_j-I'_ww_j$ and each term is zero because, as
with the $\{u_j\}$, the $\{w_j\}$ satisfy $(ww^*)_{I'}=0$ if
$|\Omega-I'|<m$. This contradiction shows that $J'=\emptyset$, and
therefore either $v_jv_i^*=0$ for all $i,j$ or $v_j=0$ for all $j$.
In the latter case, there is nothing to prove. In the former case,
since $v_j\top v_i$,
\[
v_j=v_iv_i^*v_j+v_jv_i^*v_i=(\prod_{i\in\Omega-\{j\}}v_iv_i^*)v_j
=(\prod_{i\in\Omega-\{j\}}u_iu_i^*)u_j+
(\prod_{i\in\Omega-\{j\}}w_iw_i^*)w_j=0
\]
as required.\qed

\medskip

We shall now assume that our set $\Omega$ is countable and for
convenience set $\Omega={\NN}=\{1,2,3,\ldots\}$ with its natural
order. Note that in this case, the number of possible sets $I$ in
(\ref{eq:12}) is $\aleph_0$ and the number of such $J$ is also
$\aleph_0$.

Again as in \cite{NeaRus03}, we assign a {\it signature} to each
``one'' $u_{I,k,J}$ as follows: Let the elements of $I$ be
$i_1<i_2<\cdots<i_p<\cdots$ and the elements of $J$ be
$j_1<j_2<\cdots<j_m$, where $p$ is chosen such that
$\max\{k,j_m\}<i_p$. Then $\epsilon(I,k,J)$ is defined to be the
signature of the permutation taking the $(p+m+1)$-tuple
$(i_1,\ldots,i_p,k,j_1,\ldots,j_m)$ onto $(1,2,\ldots,p+m+1)$. This
is clearly independent of $p$ as long as $\max\{k,j_m\}<i_p$.

The proof of \cite[Lemma 6.7]{NeaRus03} shows that every element
$u_{I,J}$ decomposes uniquely into a product of ``ones.'' The {\rm
signature} $\epsilon(I,J)$ (also  denoted by $\epsilon(IJ)$)
 of $u_{I,J}$ is defined to be the product of the
signatures of the factors in this  decomposition. Then the proof of
\cite[Prop.\ 6.10]{NeaRus03} shows that the family
$\{\epsilon(IJ)u_{I,J}\}$ forms a rectangular grid which satisfies
the extra property
\begin{equation}\label{eq:22}
\epsilon(IJ)u_{IJ}[\epsilon(IJ')u_{IJ'}]^*\epsilon(I'J')u_{I'J'}=
\epsilon(I'J)u_{I'J}.
\end{equation}

It follows as in \cite{NeaRus03} that the map $\epsilon(IJ)u_{IJ}
\rightarrow E_{JI}$ is a ternary isomorphism (and hence complete
isometry) from the norm closure of $\mbox{sp}_{C}\, u_{IJ}$ to the
norm closure of $\mbox{sp}_{C}\, \{E_{JI}\}$, where $E_{JI}$ denotes
an elementary matrix, whose rows and columns are indexed by the sets
$J$ and $I$, with a 1 in the $(J,I)$-position. By \cite[Lemma
1.14]{DanFri87}, this map can be extended to a ternary isomorphism
from the w*-closure of $\mbox{sp}_{C}\, u_{IJ}$ onto the Cartan
factor of type I consisting of all $\aleph_0$ by $\aleph_0$ complex
matrices which act as bounded operators on $\ell_2$. By restriction
to $Y$ and (\ref{eq:sum}), $Y$ is completely isometric to a
subtriple $\tilde Y$, of this Cartan factor of type 1.

\begin{defn}
We shall denote the space $\tilde Y$ above by $H_\infty^{m,R}$. An
entirely symmetric argument (with $J$ infinite and $I$ finite) under
the assumption that $Y$ satisfies the conditions of $Yp_L$ in Lemma
\ref{lem:3.1} part (c) with $m=m_L$  defines the space
$H_\infty^{m,L}$.
\end{defn}
Explicitly,
\[
H_\infty^{m,R}=\overline{\mbox{sp}}_{C}\, \{b_i^m=\sum_{I\cap
J=\emptyset,(I\cup
J)^c=\{i\},|J|=m}\epsilon(I,i,J)e_{J,I}:{i\in\NN}\}.
\]
and
\[
H_\infty^{m,L}=\overline{\mbox{sp}}_{C}\, \{\sum_{I\cap
J=\emptyset,(I\cup
J)^c=\{i\},|I|=m}\epsilon(I,i,J)e_{J,I}:{i\in\NN}\},
\]
with $\epsilon(I,i,J)$ defined in the obvious analogous way with $I$
finite instead of $J$.

Corollary~\ref{dist} below shows that these spaces are all distinct
from each other and from $\Phi$. This discussion has proved the
following lemma.

\begin{lemma}\label{coord}
The spaces $p_RY$ and $Yp_L$ in Lemma \ref{lem:3.1} parts (c) and
(c$^\prime$) are completely isometric to $H_\infty^{m_{R},R}$ and
$H_\infty^{m_{L},L}$, respectively.
\end{lemma}

\begin{remark}It is immediate from \cite[Cor.\ 5.3]{NeaRus03} that
$H_\infty^{0,R}=C$ and $H_\infty^{0,L}=R$.  Also note that
$H_\infty^{m,R}$ and $H_\infty^{m,L}$ are homogeneous Hilbertian
operator spaces by Lemma~\ref{lem:2.1} and \cite[Theorem
1]{NeaRus05}.
\end{remark}

\subsection{The coordinatization of Hilbertian JC*-triples}\label{sec:2.2}

Let $Y$ satisfy the hypothesis of Lemma \ref{lem:3.1}.  Our analysis
will consider the following three mutually exhaustive and (by
Corollary~\ref{dist}(b)) mutually exclusive possibilities (in each
case, the set $F$ is allowed to be empty):
\begin{description}
\item[Case 1] $(uu^*)_{\Omega-F}\ne 0$ for some finite set
$F\subset\Omega$;
\item[Case 2]
$(u^*u)_{\Omega-F}\ne 0$ for some finite set $F\subset\Omega$;
\item[Case 3]$(uu^*)_{\Omega-F}=(u^*u)_{\Omega-F}=0$ for all finite subsets
$F$ of $\Omega$.
\end{description}

We will first address cases 1 and 2.

\begin{proposition}\label{lem:3.2+3.3}
Let $Y$ be a separable infinite dimensional Hilbertian operator
space which is a $JC^*$-subtriple of $B(H)$ and let
$\{u_i:i\in\Omega\}$ ($\Omega=\NN$) be an orthonormal basis
consisting of a maximal family of mutually collinear partial
isometries of $Y$.
\begin{description}
\item[(a)]
 Suppose
there exists a finite subset $F$ of $\Omega$ such that
$(uu^*)_{\Omega-F}\ne 0$.  Then $Y$ is completely isometric to an
intersection $Y_1\cap Y_2$ such that $Y_1$ is completely isometric
to a space $H_\infty^{m,R}$  (that is, in the notation of
Lemma~\ref{lem:3.1}, $m_R^\prime=m\ge 1$ and $k_R^\prime=m+1$), or
$C$, and $Y_2$ is a Hilbertian $JC^*$-triple.
\item[(b)] Suppose  there exists a finite
subset $F$ such that $(u^*u)_{\Omega-F}\ne 0$.  Then $Y$ is
completely isometric to an intersection $Y_1\cap Y_2$ such that
$Y_1$ is completely isometric to a space $H_\infty^{m,L}$ (that is,
$m_L^\prime=m\ge 1$ and $k_L^\prime=m+1$), or $R$, and $Y_2$ is an
Hilbertian $JC^*$-triple.
\end{description}
\end{proposition}
\pf\  (a) follows from Lemma~\ref{lem:3.1} (b) and (c), the
coordinatization procedure outlined in subsection 2.1, and Lemma
\ref{coord}. (b) follows by symmetry using Lemma~\ref{lem:3.1}
(b$^\prime$) and (c$^\prime$), and Lemma \ref{coord}.\qed

\medskip


It is worth emphasizing that the space $H_\infty^{m,R}$ (resp.
$H_\infty^{m,L}$) is determined up to complete isometry among
Hilbertian JC*-triples by the condition
\begin{equation}\label{eq:m+1} (uu^*)_{\Omega-F}=0\Leftrightarrow
|F|<m,\quad (u^*u)_G=0\Leftrightarrow |G|>m+1 \end{equation}
  (resp.
$$
(u^*u)_{\Omega-F}=0\Leftrightarrow |F|<m,\quad
(uu^*)_G=0\Leftrightarrow |G|>m+1),
$$
and that $H_\infty^{0,R}=C$ and  $H_\infty^{0,L}=R$.

\begin{remark}\label{rem:m+1}
Recall from \cite[p.\ 2245]{NeaRus03} that $i_R$ (resp. $i_L$) is
the largest $i$ such that $(uu^*)_J\ne 0$ (resp.\ $(u^*u)J\ne 0$)
for any $J$ with $|J|=i$. For the spaces $H_n^k$ from
\cite{NeaRus03}, we have $i_R=k$ and $i_L=n-k+1$ so that
$i_R+i_L=n+1$.  We may therefore think of the condition
(\ref{eq:m+1}) as ``$i_R=\infty-m$, $i_L=m+1$'', so that
``$i_R+i_L=\infty+1$''.
\end{remark}

 To handle the remaining case 3, we shall need the following lemma.

\begin{lemma}\label{lem:6.0}
Let $Y$ be a separable infinite dimensional Hilbertian operator
space which is a $JC^*$-subtriple of a \csa\ $A$ and let
$\{u_i:i\in\Omega\}$ be an orthonormal basis consisting of a maximal
family of mutually collinear partial isometries of $Y$.

Let $S$ and $T$ be finite subsets of $\Omega$ and let
$k\in\Omega-(S\cup T)$. If $(uu^*)_Su_k(u^*u)_T=0$, then
$(uu^*)_{S'}u_{k'}(u^*u)_{T'}=0$ for all sets $S',T'$ with
$|S'|=|S|,\ |T'|=|T|$ and for all $k'\in \Omega-(S'\cup T')$.
\end{lemma}
\pf\ It suffices to prove this with $(S,k,T)$ replaced in turn by
$(S\cup\{l\}-\{j\},k,T)$ (with $l\not\in S$ and $j\in S$); by
$(S,l,T)$ (with $l\ne k$); and by $(S,k,T\cup\{l\}-\{i\})$ (with
$l\not\in T$ and $i\in T$).

In the first case,
\begin{eqnarray*}
u_lu_l^*(uu^*)_{S-j}u_k(u^*u)_T&=&(u_ju_j^*u_l+
u_lu_j^*u_j)u_l^*(uu^*)_{S-j}u_k(u^*u)_T\\
&=&0+u_lu_j^*u_ju_l^*(uu^*)_{S-j}u_k(u^*u)_T\\
&=&u_lu_j^*(uu^*)_{S-j}u_ju_l^*u_k(u^*u)_T\ (\mbox{by hopping})\\
&=&-u_lu_j^*(uu^*)_{S-j}u_ku_l^*u_j(u^*u)_T\\
&=&-u_lu_j^*(uu^*)_{S-j}u_k(u^*u)_Tu_l^*u_j=0.
\end{eqnarray*}
By symmetry, $(uu^*)_Su_k(u^*u)_{T-i}u_l^*u_l=0$, proving the second
case.

Finally,
\begin{eqnarray*}
(uu^*)_Su_l(u^*u)_T&=&(uu^*)_S(u_lu_k^*u_k+u_ku_k^*u_l)(u^*u)_T\\
&=&(uu^*)_Su_lu_k^*u_k(u^*u)_T+
(uu^*)_Su_ku_k^*u_l(u^*u)_T\\
&=&u_lu_k^*(uu^*)_Su_k(u^*u)_T+ (uu^*)_Su_k(u^*u)_Tu_k^*u_l=0. \qed
\end{eqnarray*}

We can now handle the final case 3.

\begin{proposition}\label{lem:3.4}
Let $Y$ be a separable infinite dimensional Hilbertian JC*-triple
and let $\{u_i:i\in\Omega\}$ be an orthonormal basis consisting of a
maximal family of mutually collinear partial isometries of $Y$.
 Suppose that for all finite subsets $G\subset\Omega$,
$(uu^*)_{\Omega-G}=0$ and $(u^*u)_{\Omega-G}=0$.  Then $Y$ is
completely isometric to $\Phi$.
\end{proposition}
\pf\ We show first that all finite products $(uu^*)_Fu_i(u^*u)_G$
with $F,G,\{i\}$ pairwise disjoint (and $F,G$ not both empty), are
not
zero. 
Suppose, on the contrary, that $(uu^*)_Fu_i(u^*u)_G=0$ for some
$F,G,i$. If $F$ and $G$ are both non-empty, pick a subset $F'\subset
F$ of maximal cardinality such that $(uu^*)_{F'}u_i(u^*u)_G\ne 0$
($F'$ could be empty). Then by repeated use of collinearity and
passing to the limit, we arrive at
$(uu^*)_{F'}u_i(u^*u)_G=(uu^*)_{F'}u_i(u^*u)_{\Omega-(\{i\}\cup
F')}=0$, a contradiction.  So either $F=\emptyset$ and $u_i(u^*u)_G=
0$, or $G=\emptyset$ and $(u^*u)_Fu_i=0$. In the first case, picking
a subset $G'\subset G$ of maximal cardinality such that
$u_i(u^*u)_{G'}\ne 0$, then by collinearity
$u_i(u^*u)_{G'}=(uu^*)_{\Omega-(\{i\}\cup G')}u_i(u^*u)_{G'}=0$, a
contradiction, and similarly in the second case. 
%
%
 We have now shown that all finite products
$(uu^*)_Fu_i(u^*u)_G$ with $F,G,\{i\}$ pairwise disjoint, are not
zero.

Now consider the space $Y_n:=\mbox{sp}\{u_1,\ldots,u_n\}$. By
\cite[Th.\ 3(b)]{NeaRus03}, $Y_n$ is completely isometric to a space
$H_n^{k_1}\cap\cdots\cap H_n^{k_m}$, where $n\ge k_1>\cdots>k_m\ge
1$. We claim that $m=n$ and $k_j=n-j+1$ for $j=1,\ldots,n$. By way
of contradiction, suppose that there is a $k$, $1\le k\le n$ such
that the space $H_n^k$ is not among the spaces $H_n^{k_j}$, $1\le
j\le m$.  Let $\psi: x\mapsto (x^{(k_1)},\dots,x^{(k_m)})$ denote
the  ternary isomorphism of the ternary envelope of $Y_n$ whose
restriction to $Y_n$ implements the complete isometry of $Y_n$ with
$H_n^{k_1}\cap\cdots\cap H_n^{k_m}$, and consider the element
$x:=(uu^*)_{\{1,\ldots,k-1\}}u_k(u^*u)_{\{k+1,\ldots,n\}}$. As shown
above, $x\ne 0$. However, $x^{(k_j)}=0$ for each $j$, a
contradiction. To see that $x^{(k_j)}=0$, suppose first that
$k_j<k$. Since $\psi$ is a ternary isomorphism,
$x^{(k_j)}=(u^{(k_j)}u^{(k_j)*})_{\{1,\ldots,k-1\}}u_k^{(k_j)}(u^{(k_j)*}u^{(k_j)})_{\{k+1,
\ldots,n\}}=0$ since
$(u^{(k_j)}u^{(k_j)*})_{\{1,\ldots,k-1\}}u_k^{(k_j)}$ is zero in
$H_n^{k_j}$. Similarly, if $k_j>k$, then $n-k_j+1<n-k+1$,
$u_k^{(k_j)}(u^{(k_j)*}u^{(k_j)})_{\{k+1, \ldots,n\}}=0$ so that
$x^{(k_j)}=0$ in this case as well.

We now have for each $n$ that, completely isometrically,
$Y_n=\cap_{k=1}^n H_n^k$ and the latter space is completely
isometric to $\Phi_n$ by \cite[Lemma 2.1]{NeaRus05}. Since
$Y=\overline{\cup Y_n}$ and $\Phi=\overline{\cup \Phi_n}$, it
follows that $Y=\Phi$ completely isometrically.
 \qed

\medskip

We come now to the first main result of this paper.
\begin{theorem}\label{thm:1}
Let $Y$ be a separable infinite dimensional Hilbertian operator
space which is concretely represented as a $JC^*$-triple. Then $Y$
is completely isometric to one of the following spaces, where each
$Y_i$ is completely isometric to one of the spaces
$H_\infty^{m_i,R}$ or $H_\infty^{m_i,L}$.
\begin{itemize}
\item $\Phi$.
\item  $Z\cap Y_1\cap\cdots \cap Y_n$, where $n\ge 1$ and $Z$ is
$\Phi$ or is absent.
\item $Z\cap Y_1\cap\cdots \cap Y_n\cap\cdots$, where $Z$ is
$\Phi$ or is absent.
\end{itemize}
\end{theorem}
\pf\ Let $\{u_i\}$ be an orthonormal basis for $Y$ consisting of a
maximal family of mutually collinear minimal partial isometries. By
Proposition~\ref{lem:3.4} and Proposition~\ref{lem:3.2+3.3}, either
$Y$ is completely isometric to $\Phi$, in which case the theorem is
proved, or $Y$ is completely isometric to an intersection $Y_1\cap
Z$, where $Y_1$ is completely isometric to either $H_\infty^{m,R}$
or $H_\infty^{k,L}$, with $m,k\ge 0$. It follows by induction, in
the case that $Y$ is not completely isometric to $ \Phi$, that
$Y=Y_1\cap Y_2\cap\cdots\cap Y_n\cap Z_n$, where $n\ge 1$, each
$Y_i$ is completely isometric to one of the spaces
$H_\infty^{m_i,R}$ or $H_\infty^{k_i,L}$, and either $Z_n$ is
completely isometric to the space $\Phi$, in which case the theorem
is proved, or $Z_n$ can be further decomposed as $Y_{n+1}\cap W$,
where $Y_{n+1}$ is completely isometric to either
$H_\infty^{m_{n+1},R}$ or $H_\infty^{k_{n+1},L}$.

It remains to consider the case in which no $Z_n$ is completely
isometric to $\Phi$. By the constructions in the proof of
Proposition~\ref{lem:3.2+3.3}, $Z_n$ is obtained by multiplying $Y$
on the left and right by a sequence of projections of the form
$1-p_{i,R}$ or $1-p_{j,L}$ (see Lemma~\ref{lem:3.1}). The resulting
products of projections converge strongly and it follows that $Y$ is
the intersection of an infinite sequence $Y_i$ and a space $Z$,
which has the property that all products $(ww^*)_{\Omega-F}$ and
$(w^*w)_{\Omega-F}$ are zero, for all finite sets $F$. An appeal to
Proposition~\ref{lem:3.4} now shows that $Z$ is completely isometric
to $\Phi$, completing the proof. \qed

\begin{remark}
By the argument in \cite[p.\ 2259]{NeaRus03}, the sequences $m_i$
and $k_i$ are strictly increasing. This fact is not needed in the
preceding proof.
\end{remark}

\begin{remark}In section 5, we show that the spaces $H_\infty^{m,R}$ and
$H_\infty^{k,L}$ are completely isometric to spaces of creation and
annihilation operators on pieces of the anti-symmetric Fock space.
Hence all rank 1 JC*-triples are really spaces of creation and
annihilation operators.
\end{remark}

We close this section with a well known lemma about Hilbertian
TRO's. Recall that TROs are operator subspaces of a C*-algebra which
are closed under the product $xy^{\ast}z$, and are fundamental in
operator space theory. 
Indeed, every operator space has both a canonical injective envelope
\cite{Ru2} and a canonical ``Shilov boundary'' \cite{Ble} which are
TROs. A proof of the following lemma can be found in \cite{Ru},
which classifies all W*-TRO's up to complete isometry. We include a
quick alternate proof from the point of view of this section.

\begin{lemma}\label{TRO}If $X$ is a Hilbertian TRO, then $X$ is
completely isometric to $R$ or $C$
\end{lemma}
\pf\ Let $\{u_j\}$ be an orthonormal basis consisting of mutually
collinear minimal partial isometries in $X$.  For a fixed $i\ne j$,
since $u_iu_i^*u_j$ is a partial isometry in $X$,
$u_iu_i^*u_j=P_2(u_j)(u_iu_i^*u_j)$ is either equal to
$e^{i\theta}u_j$ or 0.
If the latter case holds, then by the calculation in \cite[Lemma
5.1]{NeaRus03}, $u_{i}u_{i}^{\ast}u_{j}=0$ for all $i\ne j$, and $X$
is ternary isomorphic and thus completely isometric to $C$. On the
other hand, if $u_{i}u_{i}^{\ast}u_{j}=e^{i\theta}u_{j}$, then by
collinearity, $e^{i\theta}=1$, $u_{j}u_{i}^{\ast}u_{i}=0$, and again
by \cite[Lemma 5.1]{NeaRus03}, $X$ is completely isometric to $R$.
\qed

\section{Contractively complemented Hilbertian operator spaces}

Suppose that a Hilbert space $H$ is  complemented in a C*-algebra
$A$ via a contractive projection $P$.  Let $L$ be a contractive
linear map from $H$ into $A$ with the properties that $L(H) \perp H$
and $P(L(H))=0$.  Then the space $K=\{h+L(h):h \in H \}$ is clearly
contractively complemented by $P+LP$.   From this it follows that a
classification of contractively complemented Hilbertian operator
spaces is hopeless without some qualifications.

\subsection{Expansions of contractive projections}

The following definitions are crucial.

\begin{defn}Consider a triple $\{K,A,P \}$ consisting of a Hilbertian
operator space $K$, a C*-algebra $A$, and a contractive projection
$P$ from $A$ onto $K$.  If there exists a Hilbertian subspace $H$ of
$A$ which is contractively complemented by a projection $Q$ and a
contractive linear map $L$ from $H$ into $A$ such that  $P=Q+LQ$,
$L(H) \perp H$ and $Q(L(H))=0$, we say that$\{K,A,P \}$ is an {\bf
expansion} of $\{H,A,Q\}$. (Note that this implies that
$K=\{h+L(h):h \in H \}$.)
\end{defn}

The following is immediate.

\begin{lemma}
If $\{K,A,P \}$ is an expansion of $\{H,A,Q\}$ then $Q|_{K}$ is a
completely contractive isometry from $K$ onto $H$
\end{lemma}

Suppose $X \subset A$ is a contractively complemented Hilbertian
operator subspace by a projection $Q$. Further suppose that $Y$ is a
Hilbertian operator subspace of $A$ which is isometric to $X$ and
which is orthogonal to $X$ and lies in $\ker(Q)$. Then $\{x+Lx:x \in
X\}$ is contractively complemented in $A$ by the projection
$P=Q+LQ$, where $L$ is any isometry from $X$ onto $Y$. It is clear
that $\{x+Lx:x \in X\}$ is an expansion of $X$.  Thus one cannot
hope to classify contractively complemented Hilbertian operator
spaces up to complete isometry. However, we will show in this
section that all contractively complemented Hilbert spaces are
expansions of a``minimal'' 1-complemented Hilbert space which is a
JC*-triple.


\begin{defn}The {\bf support partial isometry} of a non-zero element
$\psi$ of the predual $A_{\ast}$ of a JW*-triple $A$ is the smallest
element of the set of partial isometries $v$ such that
$\psi(v)=\|\psi\|$, and is denoted by $v_\psi$.  For each non-empty
subset $G$ of $A_*$, the {\bf support space} $s(G)$ of $G$ is the
smallest weak*-closed subspace of $A$ containing the support partial
isometries of all elements of $G$.
 \end{defn}

 The existence and uniqueness of the support
partial isometry was proved and exploited in the more general case
of a \jbst\ (in which case the partial isometries are replaced by
their abstract analog, the tripotents) in \cite{FriRus85bis}. One of
its important properties is that of ``faithfulness'': if a non-zero
partial isometry $w$ satisfies $w\le v_\psi$, then $\psi(w)>0$.

We now give two examples of expansions which naturally occur and are
relevant to our work.

\begin{example}
 From \cite[Theorem 2]{NeaRus05}, if $P$ is a contractive
projection on a \csa\ $A$, with $X:=P(A)$ which is isometric to a
Hilbert space, then there are projections $p,q\in A^{**}$, such
that, $X=P^{\ast\ast}A^{\ast\ast} =\{pxq + (1-p)x(1-q):x \in X\}$.
The space $pXq$ is exactly the norm closed span  of the support
partial isometries of the elements of $P^{\ast}A^{\ast}$ (see
\cite{FriRus85} for the construction).  The map ${\bf E_0}:x\mapsto
pxq$ is an isometry of $X$ onto a $JC^*$-subtriple ${\bf E_0}X$ of
$A^{**}$, ${\bf E_0} P^{**}$ is a normal contractive projection on
$A^{**}$ with range ${\bf E_0}X$ and clearly $pXq \perp
(1-p)X(1-q)$. It follows that
$$\{X,A^{\ast\ast},P^{\ast\ast}\}\mbox{ is an expansion of }\{{\bf
E_0}X,A^{\ast\ast},{\bf E_0} P^{**}\}.$$ Specifically, let $L:{\bf
E_0}X\rightarrow A^{**}$ be the map $pxq\mapsto(1-p)x(1-q)$. Then
$P(A)=P^{**}A^{**}=\{pxq+(1-p)x(1-q):x\in P(A)\}$ and $P^{**}={\bf
E_0}P^{**}+L{\bf E_0}P^{**}$, since if $a\in A^{**}$, there is $x\in
A$ with $a=Px=p(Px)q+(1-p)Px(1-q)$ and ${\bf E_0}P^{**}a+L{\bf
E_0}P^{**}a={\bf E_0}Px+L{\bf E_0}Px=p(Px)q+(1-p)Px(1-q)$. Finally,
if $x\in A$, then $L{\bf E_0}Px=(1-p)Px(1-q)$ and ${\bf
E_0}P^{**}L{\bf E_0}Px= {\bf E_0}P^{**}((1-p)Px(1-q))=0$ by
\cite[Theorem 2(e)]{NeaRus05}.
\end{example}

\begin{defn}
The triple $\{{\bf E_0}X,A^{\ast\ast},{\bf E_0} P^{**}\}$ (or simply
${\bf E_0}X$) will be called the {\bf support} of
$\{X,A^{\ast\ast},P^{**}\}$.  It is also called the {\bf enveloping
support} of $\{X,A,P\}$.
\end{defn}

\begin{example}
It follows from \cite{EdwHugRut03}, that, for a normal contractive
projection $P$ from a von Neumann algebra (or \jwst) $A$ onto a
Hilbert space $X$, there is a similar projection ${\bf E}$ on $A$
such that $$\{X,A,P\}\mbox{ is an expansion of }\{{\bf E}A,A,{\bf E}
\}$$ and ${\bf E}A$ is the norm closure of the span of support
partial isometries of elements of $P_{\ast}A_{\ast}$.

Indeed, as set forth in \cite[Lemma 3.2]{EdwHugRut03}, $P(A)\subset
s(P_*(A_*))\oplus s(P_*(A_*))^\perp\subset A$, and ${\bf
E}:A\rightarrow A$ is a normal contractive projection onto
$s(P_*(A_*))$ given by ${\bf E}=\phi\circ P$ where
$\phi:P(A)\rightarrow s(P_*(A_*))$ is the restriction of the
$M$-projection of $s(P_*(A_*))\oplus s(P_*(A_*))^\perp$ onto
$s(P_*(A_*))$. (Although we will not use these facts, $\phi$ is a
triple isomorphism from $P(A)$ with the triple product
$\tp{x}{y}{z}_{P(A)}:=P\tp{x}{y}{z}$ onto the $JW^*$-subtriple
$s(P_*(A_*))$ of $A$, and $\phi^{-1}$ coincides with $P$ on
$s(P_*(A_*))$.

The map $L:{\bf E}(A)\rightarrow {\bf E}(A)^\perp$ in this case is
given by $L=\phi^\perp\circ \phi^{-1}=\phi^\perp\circ P$, where
$\phi^\perp:P(A)\rightarrow s(P_*(A_*))^\perp$ is the restriction of
the $M$-projection of $s(P_*(A_*))\oplus s(P_*(A_*))^\perp$ onto
$s(P_*(A_*))^\perp$. Then for $h\in s(P_*(A_*))$, say $h=\phi\circ
P(x)$ for some $x\in A$, $h+Lh=\phi(Px)+\phi^\perp(Px)=Px$ so that
$P(A)=\{h+Lh:h\in s(P_*(A_*))\}$. Furthermore, for $x\in A$, ${\bf
E}x+L{\bf E}x =\phi(Px)+\phi^\perp(Px)=Px$. It is obvious that
$L{\bf E}(A)\perp {\bf E}(A)$. Finally, for $x\in A$, ${\bf E}(L{\bf
E}x)=\phi\circ P(\phi^\perp(Px))=0$ since $Px=PPx=P\phi
Px+P\phi^\perp Px$ and $\phi\circ Px=(\phi\circ P)^2x+\phi\circ
P(\phi^\perp(Px))$.
\end{example}

\begin{defn}
By analogy with Example 1, we will call $\{{\bf E}A,A,{\bf E} \}$
(or simply
 ${\bf E}A$) the {\bf support}
of $\{X,A,P\}$ in this case. If $\{X,A,P \}$ is not the expansion of
any tuple other than itself, we say that $\{X,A,P \}$ is {\bf
essential} and that $X$ is {\bf essentially normally complemented}
in $A$.
\end{defn}

A concrete instance of Example 2 is the projection of $B(H)$ onto
$R$ (or $C$). It is easy to see that $R$ and $C$ are essentially
normally complemented in $B(H)$, as is $R\cap C$ in $B(H\oplus
H)$. (See the paragraph preceding Theorem 3 belowd.)

\begin{remark}\label{ex}
If $\{P(A),A,P\}$ is as in Example 1, then
$\{P^{**}(A^{**}),A^{**},P^{**}\}$ is as in Example 2, and the
enveloping support of $P$ is the same as the support of $P^{**}$,
since both ${\bf E_0}P(A)$ and ${\bf E}(A^{**})$ coincide with the
norm closed linear span of $A^{**}$ generated by $s(P^*(A^*))$.
\end{remark}

\begin{proposition}\label{see}
Suppose $X$ is Hilbertian and complemented in a von Neumann algebra
$A$ by a normal contractive projection $P$. Then $\{X,A,P \}$ is
essential if and only if it equals its support.
\end{proposition}
\pf\ Suppose $\{X,A,P \}$ equals its support and is the expansion of
$\{Y,A,Q\}$ given by a contractive map $L$. For each partial
isometry $v \in X$, $v=w+z$ where $w$ and $z$ are orthogonal partial
isometries, $w=Qv$, $QL=0$, $z=L(w)$and $P=Q+LQ$. Suppose $v$ is the
support partial isometry of $\psi \in P_{\ast}A_{\ast}$.  Then
\[
\psi(v)=\psi(Pv)=\psi((Q+LQ)(v))=\psi((Q+QLQ)(v))=\psi(Qv)=\psi(w),
\]
and hence $w=v$, $L= 0$ and $\{X,A,P \}$=$\{Y,A,Q\}$.  The converse
is immediate. \qed

\subsection{Operator space structure of 1-complemented
Hilbert spaces} As noted  at the beginning of the previous
subsection, we cannot classify 1-complemented Hilbert spaces up to
complete isometry. However, in Theorem 2 below, we are able to give
a classification up to ``trivial'' expansion.

 We assume in what follows that $P$ is a normal
contractive projection on a von Neumann algebra $A$, whose range
$Y=P(A)$ is a $JC^*$-subtriple of $A$ of rank one, and $\{u_i\}$ is
an orthonormal basis for $Y$ consisting of a maximal family of
minimal (in $Y$) collinear partial isometries. We shall assume for
convenience that $Y$ is infinite dimensional and separable. In
Theorem~\ref{thm:2} below, we shall also be able to handle the case
of a contractive projection on a \csa.

We know from Theorem~\ref{thm:1} that $Y$ is completely isometric to
an intersection of operator spaces $\tilde{Y}={\mathcal R}\cap
{\mathcal L}\cap \Phi$, where ${\mathcal R}= \cap_i H_\infty^{r_i}$
and ${\mathcal L}= \cap_k H_\infty^{l_k}$. Some of these spaces  may
be missing, and for short we have written
$H_\infty^{r_i}=H_\infty^{r_i,R}$ and
$H_\infty^{l_k}=H_\infty^{l_k,L}$.

As shown in section~\ref{sec:coor}, the weak*-ternary envelope of
$Y$ in $A$ is generated by the partial isometries $\{u_{I,J}\}$ and
is ternary isomorphic, hence completely isometric, to a Cartan
factor $M$ of type I which is generated by the matrix units
$\{E_{J,I}\}$. We may therefore assume that $P$ is defined on $M$
and has range $\tilde{Y}={\mathcal R}\cap {\mathcal L} \cap \Phi$,
which is a $JC^*$-subtriple of $M$. We shall identify $Y$ with
$\tilde{Y}={\mathcal R}\cap {\mathcal L} \cap \Phi$, and the
weak*-ternary envelope of $Y$ with $M$.

Note that by the definition of intersection, if the operator space
structures of $H_\infty^{r_i},H_\infty^{l_k},\Phi$ come from
embeddings $H_\infty^{r_i}\subset B(H^{r_i}),H_\infty^{l_k}\subset
B(H^{l_k}),\Phi\subset B(H^{\Phi})$, then
$$[\cap_i H_\infty^{r_i}]\cap [\cap_k H_\infty^{l_k}] \cap \Phi\subset
M\subset B([\oplus H^{r_i}]\oplus[\oplus H^{l_k}]\oplus H^{\Phi}).
$$

\begin{lemma}\label{lem:k}
If $u_j=\sum_i u_j^{r_i}+\sum_k u_j^{l_k}+u_j^\Phi$ is the
decomposition of $u_j$ into orthogonal partial isometries of
$[\cap_i H_\infty^{r_i}]\cap [\cap_k H_\infty^{l_k}] \cap \Phi$, and
if $P(u_j^{r_i})=0$ for some $i$ (resp. $P(u_j^{l_k})=0$ for some
$k$), then $u_j^{r_i}=0$ (resp. $u_j^{l_k}=0$).
\end{lemma}
\pf\ $u_j$ is the support partial isometry of some functional
$\psi_j\in Y_*=P_*(A_*)$. By the faithfulness of $\psi_j$ on its
support, $\psi_j(u_j^{r_i})=\psi_j(P(u_j^{r_i}))=0$ implies, since
$u_j^{r_i}\le u_j$, that $u_j^{r_i}=0$. Similarly for $u_j^{l_k}$.
\qed

\medskip

We again adopt the more compact notation $Iu_iJ$, used in the proof
of Lemma~\ref{lem:3.0}, for the ``one'' $(uu^*)_Iu_i(u^*u)_J$. We
note next that for $j\ne i$,
$\tpc{u_j}{Iu_iJ}{u_j}=u_j(Iu_iJ)^*u_j=u_jJu_i^*Iu_j=(J\cup\{j\}u_i^*(I\cup\{j\})=0$
since either $j\not\in I$ or $j\not\in J$. By a conditional
expectation formula in (\ref{eq:CE}),
$$0=P(\tpc{u_j}{Iu_iJ}{u_j})=\tpc{u_j}{P(Iu_iJ)}{u_j}.$$  Since
every element of $Y$ is in the closed linear span of the $u_j$, we
may write $P(Iu_iJ)=\sum_k\lambda_k^{i,J}u_k$ and thus
$0=\sum_k\overline{\lambda_k^{i,J}}\tp{u_j}{u_k}{u_j}=
\overline{\lambda_j^{i,J}}u_j$.  We conclude that
$\lambda_j^{i,J}=0$ for $j\ne i$ and hence
$P(Iu_iJ)=\lambda_{i,J}u_i$ for each ``one'' $Iu_iJ$, where we have
written $\lambda_{i,J}$ for $\lambda^{i,J}_i$.

Now suppose that $i$ is fixed and  $k\ne i$, say $k\in J$.  Then
\begin{eqnarray*}
2\tpc{u_k}{u_i}{Iu_iJ}&=&u_ku_i^*Iu_iJ+Iu_iJu_i^*u_k=u_ku_i^*Iu_iJ\mbox{
(as $k\in J$)}\\
&=&Iu_ku_i^*u_iJ \mbox{ (by ``hopping'')}\\
&=&Iu_k((J-\{k\})\cup\{i\}).
 \end{eqnarray*}

 Thus by another conditional expectation formula in (\ref{eq:CE}),
 \begin{eqnarray*}
\lambda_{k,(J-\{k\})\cup\{i\}}u_k&=&P(Iu_k((J-\{k\})\cup\{i\}))\\
&=&P(2\tpc{u_k}{u_i}{Iu_iL})\\
&=&2\tpc{u_k}{u_i}{P(Iu_iJ)}\\
&=&2\lambda_{i,J}\tp{u_k}{u_i}{u_i} =\lambda_{i,J}u_k.
\end{eqnarray*}

 Thus $\lambda_{i,J}=\lambda_{k,(J-\{k\})\cup\{i\}}$
and so $\lambda_{i,J}=\lambda$ is independent of  $i,J$ such that
$i\not\in J$ and $|J|=m$. Similarly, it can be shown that
$\lambda_{i,J}=\lambda_{k,J}$ for any $i\not\in J,k\not\in J$.

We have now shown that there is a complex number $\lambda=\lambda_m$
such that $P(Iu_jJ)=\lambda u_j$, for all partitions $I\cup\{j\}\cup
J$ of $\Omega$ with $|J|=m$.

Now, since
$P(\sum_{|J|=m}Iu_iJ)=\sum_{|J|=m}P(Iu_iJ)=\sum_{|J|=m}\lambda_m
u_i$, we must have $\lambda_m=0$ and $P(\sum_{|J|=m}Iu_iJ))=0$
unless $r_i=0$. Thus $P(u_i^{r_i})=0$ unless $m=0$. Similarly,
$P(u_i^{l_k})=0$ unless $l_k=0$. By Lemma~\ref{lem:k},
$u_j^{r_i}=u_j^{l_k}=0$ for $r_i\ne 0$, and $l_k\ne 0$. By
Theorem~\ref{thm:1}, $P(A)$ is an intersection of at most the three
spaces $R,C,\Phi$. Together with Examples 1 and 2, Remark \ref{ex},
and Proposition \ref{see} in subsection 3.1, this proves the second
main theorem of this paper.

\begin{theorem}\label{thm:2}
Suppose $Y$ is a separable infinite dimensional Hilbertian operator
space which is contractively complemented (resp. normally
contractively complemented) in a \csa\ $A$ (resp. W*-algebra $A$) by
a projection $P$.  Then, \begin{description} \item[(a)]
$\{Y,A^{\ast\ast},P^{\ast\ast}\}$ (resp. $\{Y,A,P\}$) is an
expansion of its support $\{H,A^{\ast\ast},Q\}$ (resp.
$\{H,A^{\ast\ast},Q\}$, which is essential)
\item[(b)] $H$ is contractively complemented in $A^{\ast\ast}$ (resp. $A$)
by $Q$ and is completely isometric to either $R$, $C$, $R\cap C$, or
$\Phi$.
\end{description}
\end{theorem}

This theorem says that, in $A^{\ast\ast}$, $Y$ is the diagonal of a
contractively complemented space $H$ which is completely isometric
to $R,C$, $R \cap C$ or $\Phi$ and an orthogonal degenerate space
$K$ which is in the kernel of $P$. As pointed out at the beginning
of section 3.1, this is the best possible classification.

By \cite{Y}, the range $Y$ of a completely contractive projection on
a C*-algebra is a TRO. By Lemma \ref{TRO} it follows that, if $Y$ is
Hilbertian, $Y$ is completely isometric to $R$ or $C$. This gives an
alternate proof of the result of Robertson \cite{R}, stated here for
completely contractive projections on a \csa.

Although Theorem 2 is only a classification modulo expansions, the
following Lemma shows that it is the correct analogue for
contractively complemented Hilbert spaces.

\begin{lemma}
Suppose that $\{Y,A,P\}$ is an expansion of $\{H,A,Q\}$ and that $P$
is a completely contractive projection. Then $Y$ is completely
isometric to $H$.
\end{lemma}
\pf\ By definition of expansion, in $A^{\ast\ast}$, $Y$ coincides
with $\{h+L(h):h \in H \}$, $Q+LQ=P$, $L(H) \perp H$ and
$Q(L(H))=0$. Thus, $P^{\ast\ast}|_{H}$ is a complete contraction
from $H$ onto $Y$ with completely contractive inverse $Q|_{Y}$.
Hence, $Y$ is completely isometric to $H$. \qed

\subsection{An essential contractive projection onto $\Phi$}

As noted earlier, the spaces $R$, $C$ and $R\cap C$ are each
essentially normally contractively complemented in a von Neumann
algebra.  We now proceed to show that the same holds for $\Phi$.

We begin by taking a closer look at the contractive projection
$P=P_n^k$ of the ternary envelope $T=T(H_n^k)=M_{p_k,q_k}(\CC)$ of
$H_n^k$, onto $H_n^k$. This projection and the space $H_n^k$ were
first constructed in \cite{AraFri78} and rediscovered in
\cite{NeaRus03}. By \cite[Cor.\ 7.3]{NeaRus03},
\[
P_n^kx=\frac{1}{{n-1\choose k-1}}\sum_{i=1}^n\tr(xu_i^*)u_i.
\]
Consistent with the identification of $Y$ with $\tilde{Y}$ in the
previous subsection, we let $u_i$ denote the image of the
orthonormal basis $u_i$, of a finite dimensional $JC^*$-triple of
rank one.  Thus, $u_i=\sum\epsilon(I,i,J)E_{J,I}$.

\begin{lemma}\label{lem:3.3}
The action of $P=P_n^k$ is as follows: if $x\in T$ is not a ``one'',
then $Px=0$. If $x=\epsilon(I,i,J)E_{J,I}$ is a ``one'', then
$P(\epsilon(I,i,J)E_{J,I})=\frac{1}{{n-1\choose k-1}}u_i.$
\end{lemma}
\pf\ Suppose first that $x=\epsilon(I,J)E_{J,I}\in T$ is not a
``one'', that is, $I\cap J\ne\emptyset$. Then
\begin{eqnarray*}
xu_i^*&=&\epsilon(I,J)E_{J,I}\sum_{I',J'}\epsilon(I',i,J')E_{I',J'}\\
&=&\epsilon(I,J)\sum_{I',J'}\epsilon(I',i,J')E_{J,I}E_{I',J'}\\
&=&\epsilon(I,J)\sum_{J'}\epsilon(I,i,J')E_{J,J'}.
\end{eqnarray*}
Since $J'\cap I=\emptyset$ and $J\cap I\ne\emptyset$, $J'$ is never
equal to $J$ and so $\tr(xu_i^*)=0$ and $Px=0$.

 Suppose now that $x=\epsilon(I,i,J)E_{J,I}\in T$
is a ``one''. Then for  $1\le j\le n$, $u_j=\sum_{I'\cap
J'=\emptyset}\epsilon(I',j,J')E_{J',I'}$, and as above
$xu_j^*=\epsilon(I,i,J)\sum_{J'}\epsilon(I,i,J')E_{J,J'}$. Thus,
$\tr(xu_j^*)=1$ if $j=i$ and $\tr(xu_j^*)=0$ if $j\ne i$. It follows
that $$P(\epsilon(I,i,J)E_{J,I})=\frac{1}{{n-1\choose
k-1}}\sum_j\tr(xu_j^*)u_j=\frac{1}{{n-1\choose k-1}}u_i.\qed
$$

We  proceed to construct a contractive projection defined on a TRO
$A$ which has range $\Phi$.  Since every TRO is the corner of a
\csa, we will have constructed a projection on a \csa\ with range
$\Phi$.  Now, let $u_i$ be an orthonormal basis for the Hilbertian
operator space $\Phi$ and let $H_n=\mbox{sp}\{u_1,\ldots,u_n\}$. As
noted in the proof of Proposition~\ref{lem:3.4}, $H_n=\Phi_n$ is
completely isometric to the intersection $\cap_{i=1}^n H_n^i\subset
\oplus_{k=1}^n T(H_n^i)=$ the ternary envelope $T(H_n)$ of $H_n$ in
$A$.

We construct a contractive projection $P^n$ on $T(H_n)$ with range
$H_n$ as follows. For $x=\oplus_{i=1}^n x_i\in T(H_n)$, write
$x=\sum_{i=1}^n (0\oplus\cdots\oplus x_i\oplus\cdots\oplus 0)$,
($x_i$ is in the $i^{\mbox{th}}$-position). Then define
$$P^n(x)=\sum_{i=1}^nP^n(0\oplus\cdots\oplus x_i\oplus\cdots\oplus
0):=\frac{1}{n}\sum_{i=1}^n (P_n^i(x_i),\ldots,P_n^i(x_i)).$$ Note
that since $(P_n^i(x_i),\ldots,P_n^i(x_i))$ belongs to
$H_n=\cap_{i=1}^n H_n^i$, we  shall sometimes write it as
$((P_n^i(x_i))^1,\ldots,(P_n^i(x_i))^n)$ and view $(P_n^i(x_i))^j$
as an element of $H_n^j$.

With $u_k=(u_k,\ldots,u_k)=(u_k^1,\ldots,u_k^n)=\sum_i
(0,\ldots,u_k^i,\ldots,0)$, we have
\begin{eqnarray*}
P^n(u_k)&=&\sum_i P^n((0,\ldots,u_k^i,\ldots,0))\\
&=&\frac{1}{n}\sum_i (P_n^i(u_k^i),\ldots,P_n^i(u_k^i))\\
&=&\sum_i(u_k^i,\ldots,u_k^i)/n\\
&=& \sum_i(u_k^1,\dots,u_k^n)/n\\
&=&(u_k^1,\ldots,u_k^n)=u_k.
\end{eqnarray*}

By Lemma~\ref{lem:3.3}, $P^n$ is zero an any non-``one'', so the
range of $P^n$ is $H_n$. To calculate the action of $P^n$ on
``ones'', let $x=Iu_kJ$ be such and write $x=\oplus x_i=\oplus
I^iu_k^iJ^i$, where $x_i\in H_n^i$. We claim that
\begin{equation}\label{eq:claim} P^n(x)=\frac{u_k}{n{n-1\choose
i-1}},
\end{equation} where $|I|=i-1$. Let us illustrate this first in a
specific example: Let $n=3$, $x=u_2u_2^*u_1u_3^*u_3=x_1\oplus
x_2\oplus x_3\in H_3^1\cap H_3^2\cap H_3^^3$, so that $x_1=0$,
$x_3=0$, and $i=2$. By Lemma~\ref{lem:3.3} again,
\begin{eqnarray*}
P^3(x)&=& P^3(x_1\oplus 0\oplus 0)+P^3(0\oplus x_2\oplus
0)+P^3(0\oplus
0\oplus x_3)\\
&=& \frac{1}{3}[(P_1^3(x_1),P_1^3(x_1),P_1^3(x_1))+
(P_2^3(x_2),P_2^3(x_2),P_2^3(x_2))\\
&+&(P_3^3(x_3),P_3^3(x_3),P_3^3(x_3))]\\
&=&\frac{1}{3}[(0,0,0)+(\frac{1}{2}u_1^2,\frac{1}{2}u_1^2,\frac{1}{2}u_1^2)+(0,0,0)]\\
&=&\frac{1}{3}\frac{1}{2}(u_1^2,u_1^2,u_1^2)=\frac{1}{6}u_1.
\end{eqnarray*}

In general, for $x=\oplus x_i$ as above,  $$P^n(x)=(1/n)[\sum
(P_n^i(x_i),P_n^i(x_i),P_n^i(x_i))]=(1/n)[\frac{1}{{n-1\choose
i-1}}(u_k1,\ldots,u_k^n)],$$ as required for (\ref{eq:claim}).

\begin{lemma}\label{lem:10.1}
Under the embedding $T(H_n)=\oplus_{i=1}^nT(H_n^i)\subset
\oplus_{i=1}^{n+1}T(H_n^i)=T(H_{n+1})$, given by
$x_1\oplus\cdots\oplus x_n\mapsto x_1\oplus\cdots\oplus x_n\oplus
0$, we have $P^{n+1}|T(H_n)=P^n$.
\end{lemma}
\pf\ This is obviously true for generators $u_{I,J}$ of $T(H_n)$
which are not ``ones'' since all of the $P_n^k$ and $P_{n+1}^k$
vanish on them. On the other hand, if $Iu_kJ$ is a ``one'' in
$T(H_n)$, then by collinearity
$u_k=u_ku_{n+1}^*u_{n+1}+u_{n+1}u_{n+1}^*u_k$, and by
(\ref{eq:claim}),
\begin{eqnarray*}
P^{n+1}(Iu_kJ)&=& P^{n+1}((I\cup\{n+1\})u_kJ+Iu_k(J\cup\{n+1\}))\\
&=&\frac{1}{n+1}\frac{1}{{n\choose
i}}u_k+\frac{1}{n+1}\frac{1}{{n\choose i-1}}u_k\\
&=&\frac{1}{n}\frac{1}{{n-1\choose i-1}}u_k=P^n(Iu_kJ).\qed
\end{eqnarray*}

Lemma~\ref{lem:10.1} enables the definition of a contractive
projection $P$ on a TRO $A$ which is the norm closure in
$\oplus_{i=1}^\infty M_{p_i,q_i}(\CC)$ of $\cup_{n=1}^\infty T(H_n)$
with $P(A)=\Phi$. As noted earlier, we can assume that $A$ is a
\csa. By Example 1, $\{\Phi,A^{**},P^{**}\}$ is an expansion of
$\{{\bf E}_0\Phi,A^{**},{\bf E}_0P^{**}\}$, so ${\bf
E}_0P^{**}(A^{**})={\bf E}_0\Phi$. Thus ${\bf E}_0\Phi$ is a
normally contractively complemented $JC^*$-subtriple of $A^{**}$. By
Theorem~\ref{thm:2}, ${\bf E}_0\Phi$ is completely isometric to one
of $R,C,R\cap C,\Phi$, which we shall write as $R\cap C\cap \Phi$,
with the understanding that one or two terms in this intersection
may be missing. We claim in fact that $R$ and $C$ are both missing.

\begin{lemma}\label{lem:10.2}The support ${\bf E}_0\Phi$ of
$\{\Phi,A^{\ast\ast},P^{\ast\ast}\}$ for the above construction is
completely isometric to $\Phi$.
\end{lemma}
\pf\ Because of (\ref{eq:claim}), for any partition
$\{i_1,i_2,\ldots\}\cup\{k\}\cup \{j_1,j_2,\ldots,j_m\}$ of
$\{1,2,3,\ldots\}$,
$$P^{**}(Iu_kJ)=\lim_{n\rightarrow\infty}
P^{n+m+1}(\{i_1,\ldots,i_n\}u_kJ)=\lim_{n\rightarrow\infty}
\frac{1}{n+m+1}\frac{1}{{n+m\choose n}}u_k=0.
$$

Thus, writing ${\bf E}_0u_j={\bf E}_0u_j^C+{\bf E}_0u_j^R+{\bf
E}_0u_j^\Phi$ as in the notation of Lemma \ref{lem:k} and using
Lemma \ref{lem:3.0}, $P^{\ast\ast}{\bf E}_0(u_j^C) =P^{\ast\ast}{\bf
E}_0(u_j^R)= 0$ and thus $P^{\ast\ast}({\bf E_0}u_j) =
P^{\ast\ast}({\bf E_0}u_j^\Phi)$. Recall that ${\bf E_0}(u_{j})$ is
the support partial isometry of a norm 1 element $\psi$ in
$P^{\ast}A^{\ast}$. Since $\psi({\bf E_0}u_j)=\psi(P^{**}{\bf
E_0}u_{j})=\psi(P^{**}{\bf E_0}u_{j}^\Phi)=\psi({\bf
E_0}u_{j}^\Phi)$, it follows that ${\bf E_0}(u_{j})= {\bf
E_0}u_j^\Phi$, so that ${\bf E}_0(u_j^C) ={\bf E}_0(u_j^R)= 0$,
proving that ${\bf E}_0\Phi$ is completely isometric to $\Phi$. \qed

\medskip

Since $R, C$ and $R \cap C$ are trivially contractively complemented
in $B(H)$ as spans of finite rank operators in such a way that they
clearly equal their support spaces, this proves that each of the
spaces occurring in (b) of Theorem~\ref{thm:2} are essentially
contractively complemented.

\begin{theorem}\label{thm:converse}
The operator spaces $R,C,R \cap C$ , and  $\Phi$ are each
essentially normally contractively complemented in a von Neumann
algebra.
\end{theorem}

\section{Completely bounded Banach-Mazur
distance}\label{cbbmd}

Since all of the Hilbertian operator spaces under consideration in
this paper are homogeneous (by Lemma~\ref{lem:2.1} and \cite[Theorem
1]{NeaRus05}), the completely bounded distances can be computed by
simply computing $\|\psi\|_{\mbox{cb}}\|\psi^{-1}\|_{\mbox{cb}}$ for
any fixed unitary operator between the two Hilbert spaces,
\cite[Theorem 3.1]{Zhang97}.

\begin{theorem}
For any $m,k\ge 1$,
\begin{description}
\item[(a)] $d_{\mbox{cb}}(C,H_\infty^{m,R})=\sqrt{m+1}$.
\item[(b)] $d_{\mbox{cb}}(H_\infty^{k,L},H_\infty^{m,R})=\infty$.
\item[(c)] $d_{\mbox{cb}}(\Phi,H_\infty^{m,R})=\infty$.
\end{description}
\end{theorem}
\pf\
 We first prove (a).  Let $\{u_i\}$ (resp.
$\{v_i\}$) be any orthonormal basis for $C$ (resp.
$H_\infty^{m,R}$), and let $\psi$ be the isometry that takes $u_i$
to $v_i$. For each $n>m+1$,  let
$\tilde{H}_n^1=\mbox{sp}\,\{u_1,\ldots,u_n\}$,
$\tilde{H}_{n,R}=\mbox{sp}\,\{v_1,\ldots,v_n\}$, and
$\psi^{(n)}=\psi|\tilde{H}_n^1$. Note that for $\tilde{H}_n^1$, we
have $i_R=1$ and $i_L=n$ (see Remark~\ref{rem:m+1}) so that by
\cite[Cor.\ 5.3]{NeaRus05}, $\tilde{H}_n^1$ is completely isometric
to column space $C_n=H_n^1$.  Because of this, in what follows, we
will write $H_n^1$ for $\tilde{H}_n^1$. The space $\tilde H_{n,R}$
has $i_R=m+1<n$ and $i_L=n$, and by \cite[Th.\ 3(b)]{NeaRus03}, is
completely isometric to an intersection $H_n^{k_1}\cap\cdots\cap
H_n^{k_r}$, where $m+1=k_1>k_2>\cdots>k_r$. Now, for any $p$,
\[
\|(\psi^{(n)})_p\| =\sup_{0\ne x\in
M_p(H_n^1)}\frac{\|\psi_p(x)\|_{M_p(\tilde{H}_n)}}{\|x\|_{M_p(H_n^1)}}.
\]
Let us write $x=[x_{ij}]$ with $x_{ij}\in H_n^1$ and
$y=[y_{ij}]=\psi_p(x)$, with
$y_{ij}=\psi(x_{ij})=(y_{ij}^{k_1},\ldots,y_{ij}^{k_r})    \in
\tilde{H}_n$ where $y_{ij}^{k_l}\in H_n^{k_l}$.

Now $M_p(\tilde{H}_n)\subset M_p(H_n^{k_1})\oplus\cdots\oplus
M_p(H_n^{k_r)}$, and $M_p(H_n^1)\ni [x_{ij}]\mapsto
[y_{ij}^{k_l}]\in M_p(H_n^{k_l})$ has norm $\sqrt{k_l}$ by
\cite[Lemma 3.1]{NeaRus05}. Thus $y=[y_{ij}^{k_1}]\oplus\cdots\oplus
[y_{ij}^{k_r}]$ and
\begin{eqnarray*}\|y\|_{M_p(\tilde{H}_n)}&=&\max_{1\le l\le
r}\|[y_{ij}^{k_l }]\|_{M_p (H_n^{(k_l)})}\\ &\le& \max_{1\le l\le r}
\ \sqrt{k_l}\|x\|_{M_p (H_n^1)}\\ &=&\sqrt{m+1}\|x\|_{M_p (H_n^1)}.
\end{eqnarray*}

 Thus  $\|\psi:H_n^1\rightarrow
\tilde{H}_n\|_{\mbox{cb}}\le\sqrt{m+1}$, and  by a simple
approximation argument based on the fact that $H_\infty^{m,R}$
(resp. $C$) is the norm closure of the increasing union of the
$\tilde{H}_{n,R}$ (resp. $H_n^1$), it follows that
$\|\psi\|_{\mbox{cb}}\le\sqrt{m+1}$. Moreover, equality holds.
Indeed, by the proof of \cite[Lemma 3.1]{NeaRus05}, there exists,
for each $n\ge 1$, an element $(h_{n1}^m,\ldots,h_{nn}^m)\in
M_{1,n}(H_n^1)$, such that
$\|(h_{n1}^m,\ldots,h_{nn}^m)\|_{M_{1,n}(H_n^1)}=1$ and
$\|(\psi(h_{n1}^m),\ldots,\psi(h_{nn}^m))\|_{M_{1,n}(H_n^m)}=\sqrt{m+1}$.
Then with $x_n:=\left[\begin{array}{cc}
h_{n1}^m,\ldots,h_{nn}^m&0\\
0&0 \end{array}\right] \in M_p(C)$ and $y_n=\psi_p(x_n)$, we have
$\|x_n\|_{M_p(C)}=1$ and $\|y_n\|_{M_p(H_\infty^{m,R})}=\sqrt{m+1}$,
so that $\|\psi_p\|=\sqrt{m+1}$ and
$\|\psi\|_{\mbox{cb}}=\sqrt{m+1}$.

We next show that $\|\psi^{-1}\|_{\mbox{cb}}=1$, which will complete
the proof of (a). Let $y=[y_{ij}]\in M_p(\tilde{H}_n)$  and
$x=[x_{ij}]=(\psi^{-1})_p (y)\in M_p (H_n^1)$ so that
$x_{ij}=\psi^{-1}(y_{ij})$. Then for any $1\le l\le r$, by
\cite[Lemma 3.1]{NeaRus05} and for sufficiently large $p$,
\begin{eqnarray*}
\|\psi^{-1}:\tilde{H}_n\rightarrow
H_n^1\|_{\mbox{cb}}&=&\|(\psi^{-1})_p:M_p (\tilde{H}_n)\rightarrow
M_p(H_n^1)\|\\ &=&\sup_{y\ne
0}\frac{\|(\psi^{-1})_py\|_{M_p(H_n^1)}}{\|y\|_{M_p(\tilde{H}_n)}}\\
&\le& \sup_{y\ne 0}\frac{\sqrt {\frac{n}{n-k_l+1}}\|y\|_   {M_p
(H_n^{k_l})}}{\max_{1\le q\le r}\|[y_{ij}^{k_q}]\|_{M_p(H_n^{
k_q})}}\\ &\le& \sqrt{\frac{n}{n-k_l+1}}\le
\sqrt{\frac{n}{n-m+1}}\le 1.
\end{eqnarray*}

Again, by the proof of \cite[Lemma 3.1]{NeaRus05}, for each $n\ge
1$, there exists an element $(h_{n1}^m,\ldots,h_{nn}^m)^t\in
M_{n,1}(H_n^m)$, such that
$\|(h_{n1}^m,\ldots,h_{nn}^m)^t\|_{M_{n,1}(H_n^m)}=1$ and
$\|(\psi^{-1}(h_{n1}^m),\ldots,\psi^{-1}(h_{nn}^m))^t\|_{M_{n,1}(H_n^1)}
=\sqrt{\frac{n}{n-m+1}}$. Then with
$$y_n:=\left[\begin{array}{cc} (h_{n1}^m,\ldots,h_{nn}^m)^t&0\\
0&0 \end{array}\right] \in M_p(H_\infty^{m,R})$$ and
$x_n=(\psi_p)^{-1}(y_n)$, we have $\|y_n\|_{M_p(H_\infty^{m,R})}=1$
and $\|x_n\|_{M_p(C)}=\sqrt{\frac{n}{n-m+1}}$.  Hence
$\|\psi^{-1}\|_{\mbox{cb}}=1$ and this proves (a).

We now prove (b). Let $\{u_i\}$ (resp. $\{v_i\}$) be any orthonormal
basis for  $H_\infty^{k,L}$ (resp. $H_\infty^{m,R}$), and let $\psi$
be the isometry that takes $u_i$ to $v_i$. For each $n>\max(k+1,m)$,
let $\tilde{H}_{n,L}=\mbox{sp}\,\{u_1,\ldots,u_n\}$,
$\tilde{H}_{n,R}=\mbox{sp}\,\{v_1,\ldots,v_n\}$, and
$\psi^{(n)}=\psi|\tilde{H}_{n,L}$. Note that for $\tilde{H}_{n,L}$,
we have $i_R=n$ and $i_L=k+1$ so that by \cite[Th.\ 3(b)]{NeaRus03},
$\tilde{H}_{n,L}$ is completely isometric to an intersection
$H_n^{j_1}\cap\cdots\cap H_n^{j_s}$, where $n=j_1>j_2>\cdots>j_s$.
Similarly for $\tilde{H}_{n,R}$, we have $i_R=m+1$ and $i_L=n$ so
that by \cite[Th.\ 3(b)]{NeaRus03}, $\tilde{H}_{n,R}$ is completely
isometric to an intersection $H_n^{k_1}\cap\cdots\cap H_n^{k_r}$,
where $m+1=k_1>k_2>\cdots>k_r$.

 Now, for any $p$, with $x=[x_{ij}]=(\psi)^{-1}_p(y)$,
\begin{eqnarray*}
\|((\psi^{(n)})^{-1})_p\| &=&\sup_{0\ne y\in
M_p(\tilde{H}_{n,R})}\frac{\|\psi_p^{-1}(y)\|_{M_p(\tilde{H}_{n,L})}}{\|y\|_{M_p(\tilde{H}_{n,R})}}\\
&=&\sup_{0\ne y\in M_p(\tilde{H}_{n,R})}\frac{\max_{1\le q\le s} \
|[x_{ij}^{j_q}]\|_{M_p(H_n^{j_q})}}{\max_{1\le l\le
r}\|[y_{ij}^{k_l}]\|_{M_p(H_n^{k_l})}},
\end{eqnarray*}
which, for suitable choices of $y$, as above, is greater than
\[
\frac{\max
(\sqrt{j_1},\ldots,\sqrt{j_s})}{\max(\sqrt{k_1},\ldots,\sqrt{k_r})}=\frac{\sqrt{n}}{\sqrt{m+1}}.
\]
Thus,
$\|\psi^{-1}\|_{\mbox{cb}}\ge\|(\psi^{(n)})^{-1}\|_{\mbox{cb}}\ge\|((\psi^{(n)})^{-1})_p\|\ge
\frac{\sqrt{n}}{\sqrt{m+1}}\rightarrow\infty$. This proves (b).

Finally, we prove (c).
 Let $\{u_i\}$ (resp. $\{v_i\}$) be any
orthonormal basis for $\Phi$ (resp.\  $H_\infty^{m,R}$, and let
$\psi$ be the isometry that takes $u_i$ to $v_i$. For each $n>m$,
let $\tilde{H}_{n,R}=\mbox{sp}\,\{u_1,\ldots,u_n\}$,
$\tilde{H}_{n,\Phi}  =\mbox{sp}\,\{v_1,\ldots,v_n\}$, and
$\psi^{(n)}=\psi|\tilde{H}_{n,\Phi}$. Note that for
$\tilde{H}_{n,R}$, we have $i_R=m+1$ and $i_L=n$ so that by
\cite[Th.\ 3(b)]{NeaRus03}, $\tilde{H}_{n,R}$ is completely
isometric to an intersection $H_n^{j_1}\cap\cdots\cap H_n^{j_s}$,
where $m+1=j_1>j_2>\cdots>j_s$. For $\tilde{H}_{n,\Phi}$, we have
$i_R=n$ and $i_L=n$ so that by \cite[Th.\ 3(b)]{NeaRus03},
$\tilde{H}_{n,\Phi}$ is completely isometric to an intersection
$H_n^{k_1}\cap\cdots\cap H_n^{k_r}$, where $n=k_1>k_2>\cdots>k_r$
(in fact, as shown in the proof of Proposition~\ref{lem:3.4}, $r=n$
and $k_j=n-j+1$ but we do not need this fact).

 Now, for any $p$, with $x=[x_{ij}]\in M_p(\tilde{H}_{n,\Phi}),\ y=[y_{ij}]=\psi_p(x)$,
\begin{eqnarray*}
\|(\psi^{(n)})_p\| &=&\sup_{0\ne x\in
M_p(\tilde{H}_{n,R})}\frac{\|\psi_p(x)\|_{M_p(\tilde{H}_{n,R})}}{\|x\|_{M_p(\tilde{H}_{n,\Phi})}}\\
&=&\sup_{0\ne x\in M_p(\tilde{H}_{n,R})}\frac{\max_{1\le l\le r} \
|[y_{ij}^{k_l}]\|_{M_p(H_n^{k_l})}}{\max_{1\le q\le
s}\|[x_{ij}^{j_q}]\|_{M_p(H_n^{j_s})}},
\end{eqnarray*}
which, for suitable choices of $x$, as above, is $\ge
\frac{\max(\sqrt{k_1},\ldots,\sqrt{k_r})}{\max
(\sqrt{j_1},\ldots,\sqrt{j_s})}=\frac{\sqrt{n}}{\sqrt{m+1}}, $
 showing $\|\psi\|_{\mbox{cb}}=\infty$ and proving (c).\qed


\begin{corollary}\label{dist} For $m,k\ge 0$,
\begin{description}
\item[(a)] $H_\infty^{m,R}$ and $H_\infty^{k,R}$ are completely
isomorphic but not completely isometric if $m\ne k$.
\item[(b)] $H_\infty^{k,L}$ and $H_\infty^{m,R}$ are not
completely isomorphic.
\item[(c)] $\Phi$ and $H_\infty^{m,R}$ are not completely isomorphic.
\end{description}
\end{corollary}

Similar arguments yield the following distances as well as their
corresponding consequences: \begin{itemize}
\item $d_{\mbox{cb}}(R,H_\infty^{k,L})=\sqrt{k+1}$.
\medskip
\item $d_{\mbox{cb}}(R,H_\infty^{m,R})=d_{\mbox{cb}}(C,H_\infty^{m,L})
=d_{\mbox{cb}}(\Phi,H_\infty^{m,L})=\infty$
\end{itemize}

\section{Representation on the Fock space}

We begin by recalling the construction of the spaces
$H_\infty^{m,R}$; see subsection \ref{sec:2.1}. Let $I$ denote a
subset of $\Omega$ with $|\Omega-I|=m+1$. and let $J$ denote a
subset of $\Omega$ of cardinality $|J|=m$. We assume that each
$I=\{i_1,i_2,\ldots\}$ is such that $i_1<i_2<\cdots $, and that the
collection
 of all such subsets $I$ is ordered
lexicographically. Similarly, if $J=\{j_1,\ldots,j_{m}\}$, then
$j_1<\cdots<j_{m}$ and the collection of all such subsets $J$ is
ordered lexicographically.

We shall use the notation $e_i$ to denote the column vector with a 1
in the $i^{\mbox{th}}$ position and zeros elsewhere. Thus
$e_1,e_2,\ldots$ denotes the canonical basis of column vectors for
separable column space $C$. More generally, $e_I$ denotes the basis
vector for $\ell_2$ consisting of a 1 in the ``$I^{\mbox{th}}$''
position.

The space $H_\infty^{m,R}$ is the closed linear span of matrices
$b_i^m$, $i\in\Omega$, given by
\[
b_i^m=\sum_{I\cap J=\emptyset,(I\cup
J)^c=\{i\},|J|=m}\epsilon(I,i,J)e_{J,I},
\]
where $e_{J,I}=e_J\otimes e_I=e_Je_I^t\in
M_{\aleph_0,\aleph_0}(\CC)=B(\ell_2)$, and $\epsilon(I,i,J)$  is the
signature of the permutation defined for disjoint $I,J$ in
subsection ~\ref{sec:2.1}.

Let $H$ be any separable infinite dimensional Hilbert space. For any
$h\in H$, let $ C^m_{h}$ denote the wedge (or creation) operator
from $\wedge^{m}H$ to $\wedge^{m+1}H$ given by
$$C^m_{h}(h_1\wedge\cdots\wedge h_{m})=h\wedge
h_1\wedge\cdots\wedge h_{m}.
$$
The space of creation operators
$\overline{\mbox{sp}}\{C_{e_i}^{m}\}$ will be denoted by ${\mathcal
C}^{m}$. Its operator space structure is given by its embedding in
$B(\wedge^mH,\wedge^{m+1}H)$.

It will be convenient to identify the space $\wedge^kH$ with
$\ell_2(\{J\subset\Omega:|J|=k\})$ or with
$\ell_2(\{I\subset\Omega:|\Omega-I|=k\})$.

Define the unitary operators $V_k$ and $W_k$ on $\wedge^kH$ by
$$V_k:\ell_2(\{I\subset\Omega:|\Omega-I|=k\})\rightarrow
\ell_2(\{K\subset \Omega:|K|=k\})
$$
and
$$W_k:\ell_2(\{J\subset\Omega:|J|=k\})\rightarrow
\ell_2(\{J\subset\Omega:|J|=k\})
$$  as follows:
\begin{itemize}
\item $V_k(e_I)=e_{{\NN}-I}$; More specifically, $V_k(e_I)=e_{j_1}\wedge\cdots\wedge e_{j_k}$
where ${\NN}-I=\{j_1<\cdots<j_k\}$.
\item $W_k(e_I)=\epsilon(i,I)\epsilon(I,i,J)e_I$ where
$I\cup\{i\}\cup J=\NN$ is a disjoint union.
\end{itemize}
  It is easy to see, as in
\cite[section 2]{NeaRus05}, that the definition of $W_k$ is
independent of the choice of $i$. Indeed, if $p$ is chosen so that
$i_p>\max\{i,j_{k-1}\}$, then
$\epsilon(I,i,J)=(-1)^p\epsilon(i,i_1,\ldots,i_p,j_1,\ldots,j_{k-1})=
(-1)^p\epsilon(i,i_1,\ldots,i_p)\epsilon(\{i\}\cup\{i_1,\ldots,i_p\},j_1,\ldots,j_{k-1})$
and therefore  for any $i'\ne i$,
$\epsilon(i,I)\epsilon(I,i,J)=\epsilon(i',I')\epsilon(I',i',J)$.

\begin{lemma}\label{lem:5.1}
$H_\infty^{m,R}$ is completely isometric to ${\mathcal C}^{m}$.
\end{lemma}
\pf\
 With $b_i^m=\sum\epsilon(I,i,J)e_{J,I}$ we have
 \[
b_i^mW_m(e_{I_0})=b_i^m((\epsilon(i,I_0)\epsilon(I_0,i,J_0)e_{I_0})=\epsilon(i,I_0)e_{J_0},
 \]
and
\[
V_{m+1}C_{e_i}^m(e_{I_0})=V_{m+1}(\epsilon(i,I_0)e_{\{i\}\cup
I_0})=\epsilon(i,I_0)e_{J_0}.\qed
\]

\medskip

Since the anti-creation operator space $A^{m}$ is simply the adjoint
of the creation operator space $C^{m}$, by construction it is clear
that
\begin{lemma}\label{lem:5.2}
$H_\infty^{m,L}$ is completely isometric to the space of
anti-creation operators ${\mathcal A}^{m}$.
\end{lemma}


 By Lemma~\ref{lem:5.1}, Lemma \ref{lem:5.2}, Theorem \ref{thm:1}, and
 \cite[Lemma 2.1]{NeaRus05} we now have
\begin{theorem}\label{thm:5}Every $n$-dimensional or infinite dimensional
separable Hilbertian JC*-triple is completely isometric to an
intersection over a set of values of $k$ of the spaces of creation
and annihilation operators on $k$-fold antisymmetric tensors of the
Hilbert space.
\end{theorem}

\bibliographystyle{amsplain}

\end{document}